\documentclass[11pt]{article}

\usepackage{rotating}
\usepackage{amsmath}
\usepackage{amsfonts}
\usepackage{amssymb}
\usepackage{latexsym}

\newtheorem{thm}{Theorem}[section]

\newtheorem{lem}{Lemma} [section]

\newtheorem{cnj}[thm]{Conjecture}

\def \cB {{\cal B}}

\begin{document}

\pagestyle{headings}

\newcommand{\qed}{\hfill{\setlength{\fboxsep}{0pt}
                  \framebox[7pt]{\rule{0pt}{7pt}}} \newline}

\def \cC {{\cal C}}
\def \cD {{\cal D}}
\def \cB {{\cal B}}
\def \cE {{\cal E}}
\def \H {{\cal H}}
\def \HD {${\cal H}^*$}
\def \HH {$\H=(X,E)$}
\def \HHD {$\H^*=(E^*,X^*)$}
\def \SH {$\H=(V,{\cal B})$}
\def \cP {{\cal P}}

\date{August 22, 2013}

\title {Appendix of Extending Bicoloring for Steiner Triple Systems}

\author
{ M. Gionfriddo*\,  \and E. Guardo*\, \and L. Milazzo\,
\thanks{~ Department of Mathematics and Computer Science,
University of Catania, Viale A. Doria, 6 95125 - Catania, Italy.
E-mail: {\tt gionfriddo@dmi.unict.it}, {\tt guardo@dmi.unict.it},
{\tt milazzo@dmi.unict.it}.}}

\maketitle

\bigskip  \bigskip

\begin{center}
{\it In memory of Lucia Gionfriddo} \end{center}

\bigskip  \bigskip

\begin{abstract}
\noindent  We initiate the study of extended bicolorings of Steiner triple systems (STS) which start with a $k$-bicoloring of an STS($v$)
and end up with a $k$-bicoloring of an STS($2v+1$) obtained by a doubling construction, using only the original colors used in coloring the
subsystem STS($v$).
By producing many such extended bicolorings, we obtain several infinite classes of orders for which there exist STSs
with different lower and upper chromatic number.
\end{abstract}

\bigskip  \bigskip  \bigskip

\section{Appendix}

\noindent In this appendix we list all the factorizations $\cal F$
which characterize all the extended colorings shown in section~$3$
of the paper. In the following, until table~21, we will insert
the color classes before of the factorizations $\cal F$. All
the color classes for the tables~22, 23,24,25 and 26 are
defined at page~7.

\bigskip   \bigskip  \bigskip

\footnotesize{
\noindent \textbf{BSTS($27$)}

\noindent  Color classes:  $X_1=\{1,2,14,$ $15,16,17\}$, $X_2=\{3,4,5,$ $6,7,18,$ $19,20,21\}$ and  $X_3=\{8,9,10,$ $11,12,13,$ $22,23,24,$ $25,26,27\}$.}

\bigskip   \bigskip

\footnotesize{
\begin{center}
\begin{tabular}{ccccccc}  \label{tbl1}
1&2&3&4&5&6&7\\
14~~18&14~~21&18~~16&18~~17&18~~22&18~~23&18~~25\\
15~~19&15~~18&19~~17&19~~14&19~~23&19~~26&19~~22\\
16~~20&16~~19&14~~20&20~~15&20~~26&20~~24&20~~27\\
17~~21&17~~20&15~~21&21~~16&21~~27&21~~25&21~~24\\
22~~23&22~~24&22~~25&22~~26&14~~15&14~~16&14~~17\\
24~~27&25~~23&24~~26&25~~27&16~~17&15~~17&15~~16\\
25~~26&26~~27&23~~27&23~~24&24~~25&22~~27&23~~26\\
\end{tabular}

\bigskip  \bigskip

\begin{tabular}{cccccc}
8&9&10&11&12&13\\
22~~14&22~~15&22~~20&22~~16&22~~17&22~~21\\
23~~15&23~~16&23~~21&23~~14&23~~20&23~~17\\
24~~16&24~~17&24~~14&24~~18&24~~19&24~~15\\
25~~17&25~~14&25~~15&25~~20&25~~16&25~~19\\
26~~18&26~~21&26~~16&26~~17&26~~15&26~~14\\
27~~19&27~~18&27~~17&27~~15&27~~14&27~~16\\
20~~21&19~~20&18~~19&19~~21&18~~21&18~~20\\
~~&~~&~~&~~&~~&~~\\
\end{tabular}

table~16
\end{center}}

\bigskip

\footnotesize{
\noindent  Color classes: $X_1=\{1,2,14,$ $15,16,$ $17,18,$ $19,20\}$, $X_2=\{3,4,5,6,$ $7,21\}$ and $X_3=\{8,9,10,$ $11,12,13,$ $22,23,24,25,$ $26,27\}$.}

\bigskip

\footnotesize{
\begin{center}

\begin{tabular}{ccccccc}
1&2&3&4&5&6&7\\
14~~21&14~~23&22~~23&22~~24&22~~25&22~~26&22~~27\\
15~~27&15~~26&24~~27&23~~25&24~~26&25~~27&23~~26\\
16~~24&16~~27&25~~26&26~~27&23~~27&23~~24&24~~25\\
17~~23&17~~25&14~~16&15~~17&16~~18&17~~19&18~~20\\
18~~25&18~~24&17~~20&14~~18&15~~19&16~~20&14~~17\\
19~~26&19~~22&18~~19&19~~20&14~~20&14~~15&15~~16\\
20~~22&20~~21&21~~15&21~~16&21~~17&21~~18&21~~19\\
\end{tabular}

\bigskip

\begin{tabular}{cccccc}
8&9&10&11&12&13\\
22~~21&22~~14&22~~15&22~~16&22~~17&22~~18~~~~~~\\
23~~18&23~~21&23~~16&23~~15&23~~19&23~~20~~~~~~\\
24~~14&24~~15&24~~21&24~~17&24~~20&24~~19~~~~~~\\
25~~16&25~~20&25~~19&25~~21&25~~14&25~~15~~~~~~\\
26~~17&26~~18&26~~14&26~~20&26~~16&26~~21~~~~~~\\
27~~19&27~~17&27~~20&27~~18&27~~21&27~~14~~~~~~\\
15~~20&16~~19&17~~18&14~~19&15~~18&16~~17~~~~~~\\
~~&~~&~~&~~&~~&~~\\
\end{tabular}

table~17
\end{center}}

\footnotesize{
\noindent \textbf{BSTS($43$)}

 \noindent  Color classes: $X_1=\{1,2,3,4,5,22,23,24,25,26\}$, $X_2=\{6,7,8,$ $9,10,11,$ $27,28,29,$ $30,31,32,33,$ $34,35,$ $36\}$, $X_3=\{12,13,14,$ $15,16,17,$ $18,$ $19,$ $20,$ $21,$ $37,$ $38,39,40,41,42,43\}$.}

\bigskip
\footnotesize{
\begin{center}
\begin{tabular}{ccccccccc}  \label{tbl1}
1&2&3&4&5&6&7&8&9\\
22~~27&23~~29&26~~28&26~~36&25~~33&23~~27&24~~27&25~~27&26~~27\\
23~~28&24~~31&25~~34&25~~35&24~~34&22~~28&25~~28&24~~32&22~~32\\
25~~29&22~~35&24~~35&24~~28&23~~35&24~~29&26~~29&26~~33&23~~33\\
26~~30&26~~34&23~~36&22~~29&22~~36&25~~30&22~~30&23~~30&24~~30\\
24~~43&25~~39&22~~42&23~~41&26~~40&26~~32&23~~31&22~~31&25~~31\\
31~~34&27~~36&27~~31&27~~32&31~~32&31~~37&32~~38&28~~39&28~~38\\
32~~33&28~~32&30~~32&31~~33&28~~29&33~~41&35~~43&29~~42&29~~43\\
35~~36&30~~33&29~~33&30~~34&27~~30&34~~42&33~~40&34~~37&34~~40\\
37~~42&38~~40&41~~43&40~~42&39~~41&35~~40&34~~41&35~~38&35~~42\\
38~~41&37~~41&37~~40&39~~43&38~~42&36~~39&36~~37&36~~43&36~~41\\
39~~40&42~~43&38~~39&37~~38&37~~43&38~~43&39~~42&40~~41&37~~39\\
\end{tabular}
\end{center}}

\footnotesize{
\begin{center}
\begin{tabular}{ccccccccc}
10&11&12&13&14&15&16&17&18\\
23~~32&22~~33&22~~39&22~~37&22~~43&22~~38&22~~41&30~~37&34~~39\\
24~~33&25~~32&23~~42&23~~38&23~~39&23~~43&23~~40&31~~38&33~~42\\
22~~34&23~~34&24~~40&24~~41&24~~42&24~~39&24~~37&27~~39&32~~41\\
26~~35&24~~36&25~~41&25~~42&25~~37&25~~40&25~~38&32~~43&31~~40\\
25~~36&26~~31&26~~37&26~~43&26~~38&26~~42&26~~39&29~~41&30~~43\\
28~~37&28~~40&33~~38&31~~39&28~~41&29~~37&30~~42&28~~42&27~~38\\
29~~38&29~~39&34~~43&32~~40&29~~40&30~~41&31~~43&22~~40&23~~37\\
30~~39&30~~38&27~~29&27~~28&27~~33&27~~34&27~~35&23~~26&22~~24\\
31~~41&35~~37&28~~30&29~~36&32~~34&33~~35&34~~36&24~~25&25~~26\\
27~~42&27~~43&31~~36&30~~35&31~~35&32~~36&28~~33&33~~36&29~~35\\
40~~43&41~~42&35~~32&33~~34&30~~36&28~~31&29~~32&34~~35&28~~36\\
\end{tabular}

\begin{tabular}{ccc}
19&20&21\\
36~~42&38~~36&36~~40\\
35~~41&35~~39&34~~38\\
33~~39&33~~37&33~~43\\
32~~37&32~~42&32~~39\\
28~~43&30~~40&27~~37\\
27~~40&27~~41&26~~41\\
24~~38&25~~43&22~~25\\
23~~25&24~~26&23~~24\\
22~~26&22~~23&29~~30\\
30~~31&29~~31&28~~35\\
29~~34&28~~34&31~~42\\
~~&~~&~~\\
\end{tabular}

table~18
\end{center}}

\newpage

\footnotesize{
 \noindent  Color classes: $X_1=\{1,2,3,$ $4,22,23,$ $24,25,26\}$, $X_2=\{5,6,$ $7,8,$ $9,10,$ $11,12,$ $27,28,29,$ $30,31,32,$ $33,34,35,$ $36,37\}$, $X_3=\{13,$ $14,15,$ $16,17,$ $18,19,$ $20,21,38,$ $39,40,41,$ $42,43\}$.}

 \bigskip

\footnotesize{
\begin{center}
\begin{tabular}{ccccccccc}  \label{tbl1}
1&2&3&4&5&6&7&8&9\\
22~~27&22~~37&22~~28&27~~32&22~~29&27~~41&27~~43&27~~38&27~~25\\
23~~28&23~~30&23~~27&23~~29&23~~31&28~~26&28~~42&28~~39&28~~43\\
24~~29&24~~28&24~~31&24~~30&24~~32&29~~38&29~~40&29~~41&29~~42\\
25~~36&25~~29&25~~30&25~~28&25~~33&30~~22&30~~39&30~~40&30~~41\\
26~~37&26~~27&26~~29&26~~31&26~~30&31~~25&31~~22&31~~42&31~~38\\
38~~39&38~~40&38~~41&38~~42&38~~43&32~~23&32~~25&32~~43&32~~26\\
40~~43&39~~41&40~~42&41~~43&39~~42&33~~39&33~~41&33~~24&33~~22\\
41~~42&42~~43&39~~43&39~~40&40~~41&34~~40&34~~23&34~~22&34~~24\\
30~~35&31~~36&32~~37&27~~33&28~~34&35~~24&35~~26&35~~23&35~~39\\
31~~34&32~~35&33~~36&34~~37&27~~35&36~~42&36~~24&36~~26&36~~23\\
32~~33&33~~34&34~~35&35~~36&36~~37&37~~43&37~~38&37~~25&37~~40\\
\end{tabular}

\bigskip

\begin{tabular}{ccccccccc}
10&11&12&13&14&15&16&17&18\\
27~~42&27~~24&38~~23&38~~22&38~~26&38~~25&38~~35&38~~24&38~~36\\
28~~41&28~~38&39~~32&39~~24&39~~34&39~~23&39~~22&39~~25&39~~27\\
29~~39&29~~43&40~~24&40~~23&40~~22&40~~26&40~~28&40~~36&40~~35\\
30~~43&30~~42&41~~25&41~~26&41~~24&41~~22&41~~32&41~~23&41~~34\\
31~~40&31~~39&42~~26&42~~25&42~~23&42~~35&42~~37&42~~22&42~~33\\
32~~38&32~~40&43~~22&43~~33&43~~25&43~~24&43~~31&43~~26&43~~23\\
33~~23&33~~26&31~~33&32~~34&33~~35&34~~36&27~~36&35~~37&28~~37\\
34~~26&34~~25&30~~34&31~~35&32~~36&33~~37&29~~34&34~~27&30~~32\\
35~~25&35~~22&29~~35&30~~36&31~~37&27~~32&30~~33&28~~33&29~~31\\
36~~22&36~~41&28~~36&29~~37&27~~30&28~~31&23~~26&29~~32&22~~24\\
37~~24&37~~23&27~~37&27~~28&28~~29&29~~30&24~~25&30~~31&25~~26\\
\end{tabular}
\end{center}}

\bigskip

\footnotesize{
\begin{center}
\begin{tabular}{ccc}
19&20&21\\
38~~34&38~~33&38~~30\\
39~~36&39~~37&39~~26\\
40~~33&40~~25&40~~27\\
41~~31&41~~35&41~~37\\
42~~24&42~~32&42~~34\\
43~~35&43~~34&43~~36\\
27~~29&28~~30&29~~33\\
30~~37&27~~31&31~~32\\
28~~32&29~~36&28~~35\\
23~~25&24~~26&22~~25\\
22~~26&22~~23&23~~24\\
~~&~~&~~\\
\end{tabular}

table~19

\end{center}}

\newpage

\footnotesize{
\noindent  Color classes: $X_1=\{1,2,3,4,$ $22,23,24,$ $25,26,27\}$, $X_2=\{5,6,7,$ $8,9,10,11,$ $12,28,29,$ $30,31,32,$ $33,34,$ $35\}$, $X_3=\{13,14,15,$ $16,17,18,$ $19,20,21,$ $36,37,38,$ $39,40,41,$ $42,43\}$.}

\bigskip

\footnotesize{
\begin{center}
\begin{tabular}{ccccccccc}  \label{tbl1}
1&2&3&4&5&6&7&8&9\\
22~~28&22~~29&22~~30&22~~31&35~~36&28~~41&28~~38&32~~42&32~~39\\
23~~29&26~~28&23~~31&23~~32&29~~38&29~~36&29~~39&34~~43&29~~40\\
24~~30&24~~31&27~~29&24~~33&30~~39&30~~40&30~~37&30~~36&33~~41\\
25~~31&25~~32&25~~33&25~~34&31~~40&31~~37&31~~42&31~~38&31~~36\\
26~~40&23~~41&26~~42&26~~37&32~~22&32~~27&32~~26&28~~23&28~~25\\
27~~41&27~~42&24~~43&27~~38&33~~23&33~~22&33~~27&33~~26&35~~22\\
32~~33&30~~34&28~~35&28~~29&28~~27&35~~23&35~~24&29~~24&34~~27\\
34~~35&33~~35&32~~34&30~~35&34~~26&34~~24&34~~22&35~~25&30~~26\\
36~~43&36~~38&36~~39&36~~41&37~~43&38~~43&36~~40&37~~40&37~~38\\
37~~42&37~~39&38~~40&40~~42&41~~42&39~~42&41~~43&39~~41&42~~43\\
38~~39&40~~43&37~~41&39~~43&24~~25&25~~26&23~~25&22~~27&23~~24\\
\end{tabular}
\end{center}}

\bigskip

\footnotesize{
\begin{center}
\begin{tabular}{ccccccccc}
10&11&12&13&14&15&16&17&18\\
28~~39&32~~41&28~~42&28~~36&28~~37&28~~40&28~~43&29~~41&29~~42\\
29~~43&34~~37&34~~40&29~~37&35~~39&35~~41&33~~39&35~~37&32~~38\\
32~~37&35~~43&33~~43&35~~38&34~~38&34~~36&35~~40&33~~40&33~~37\\
33~~43&28~~24&31~~39&34~~39&32~~40&32~~43&31~~41&30~~42&30~~43\\
30~~25&29~~25&35~~27&23~~40&22~~36&27~~37&26~~36&27~~43&27~~40\\
31~~27&30~~27&29~~26&24~~42&24~~41&25~~38&22~~37&25~~39&25~~36\\
35~~26&31~~26&32~~24&26~~41&23~~42&23~~39&24~~38&24~~36&24~~39\\
34~~23&36~~42&30~~23&22~~43&25~~43&22~~42&25~~42&23~~38&22~~41\\
40~~41&39~~40&36~~37&25~~27&26~~27&24~~26&23~~27&22~~26&23~~26\\
38~~42&22~~23&38~~41&30~~33&29~~33&29~~30&29~~34&28~~32&28~~34\\
22~~24&33~~38&22~~25&31~~32&30~~31&31~~32&30~~32&31~~34&31~~35\\
\end{tabular}
\end{center}}

\bigskip

\footnotesize{
\begin{center}
\begin{tabular}{ccc}
19&20&21\\
30~~41&33~~42&35~~42\\
32~~36&34~~41&30~~38\\
34~~42&36~~23&36~~27\\
31~~43&37~~25&37~~24\\
26~~38&38~~22&39~~26\\
25~~40&39~~27&40~~22\\
23~~37&40~~24&41~~25\\
22~~39&43~~26&43~~23\\
24~~27&29~~31&28~~31\\
29~~35&28~~30&33~~34\\
28~~33&32~~35&29~~32\\
~~&~~&~~\\
\end{tabular}

table~20
\end{center}}

 \bigskip
\footnotesize{
\noindent  Color classes: $X_1=\{1,2,3,$ $4,22,23,$ $24,25,26,$  $27\}$, $X_2=\{5,6,$ $7,8,$ $9,10,$ $11,12,$ $28,29,$ $30,31,32,$ $33,34,35,$ $36\}$, $X_3=\{13,$ $14,15,$ $16,17,$ $18,19,$ $20,21,37,$ $38,39,40,$ $41,42,$. $43\}$}

\bigskip
\footnotesize{
\begin{center}
\begin{tabular}{ccccccccc}  \label{tbl1}
1&2&3&4&5&6&7&8&9\\
22~~30&27~~28&24~~35&26~~36&28~~23&36~~22&28~~22&29~~22&31~~22\\
23~~34&23~~38&22~~37&22~~40&29~~24&33~~23&30~~23&31~~23&36~~23\\
25~~32&24~~40&23~~41&23~~42&30~~25&32~~24&34~~24&30~~24&33~~24\\
24~~41&25~~39&25~~42&24~~43&31~~26&34~~25&36~~25&35~~25&29~~25\\
26~~37&26~~42&26~~43&25~~37&32~~27&28~~26&29~~26&34~~26&35~~26\\
27~~40&22~~43&27~~38&27~~39&33~~22&31~~27&35~~27&33~~27&34~~27\\
38~~43&37~~41&39~~40&38~~41&34~~37&29~~43&31~~38&28~~37&28~~41\\
39~~42&29~~36&28~~30&28~~33&35~~40&30~~37&32~~39&32~~41&30~~39\\
29~~31&30~~35&29~~33&29~~34&36~~42&35~~42&33~~43&36~~38&32~~38\\
33~~36&31~~34&31~~36&30~~32&38~~39&39~~41&37~~40&39~~43&37~~42\\
28~~35&32~~33&32~~34&31~~35&41~~43&38~~40&41~~42&40~~42&40~~43\\
\end{tabular}
\end{center}}

\bigskip

\footnotesize{
\begin{center}
\begin{tabular}{ccccccccc}
10&11&12&13&14&15&16&17&18\\
32~~22&34~~22&35~~22&37~~29&37~~31&37~~32&37~~35&37~~36&38~~22\\
35~~23&32~~23&29~~23&38~~30&38~~28&38~~33&38~~34&38~~35&37~~23\\
28~~24&36~~24&31~~24&39~~31&39~~33&39~~35&39~~29&39~~34&42~~24\\
33~~25&31~~25&28~~25&40~~32&40~~34&40~~36&40~~28&40~~33&41~~25\\
30~~26&33~~26&32~~26&41~~33&41~~36&41~~29&41~~30&41~~31&39~~26\\
36~~27&29~~27&30~~27&42~~34&42~~30&42~~28&42~~32&42~~29&43~~27\\
29~~38&28~~39&36~~39&43~~35&43~~32&43~~31&43~~36&43~~30&29~~40\\
31~~40&30~~40&33~~42&22~~23&22~~24&22~~25&22~~26&22~~27&28~~31\\
34~~41&35~~41&34~~43&24~~27&23~~25&24~~26&25~~27&23~~26&30~~36\\
37~~39&38~~42&37~~38&25~~26&26~~27&23~~27&23~~24&24~~25&32~~35\\
42~~43&37~~43&40~~41&28~~36&29~~35&30~~34&31~~33&28~~32&33~~34\\
\end{tabular}
\end{center}}

\bigskip
\footnotesize{
\begin{center}
\begin{tabular}{ccc}
19&20&21\\
42~~22&41~~22&39~~22\\
39~~23&43~~23&40~~23\\
37~~24&39~~24&38~~24\\
40~~25&38~~25&43~~25\\
38~~26&40~~26&41~~26\\
41~~27&42~~27&37~~27\\
28~~43&33~~37&31~~42\\
29~~30&28~~34&28~~29\\
31~~32&29~~32&30~~33\\
33~~35&30~~31&34~~35\\
34~~36&35~~36&32~~36\\
~~&~~&~~\\
\end{tabular}

table~21
\end{center}}

\bigskip
\newpage

\footnotesize{

\noindent \textbf{BSTS($91$)}

\noindent   Color classes of Table~22,: $X_1=\{1,2,3,$ $4,46,47,$ $48,49\}$, $X_2=\{5,6,7,$ $8,9,10,$ $50,51,52,$  $53,54,55,$ $56,57\}$, $X_3=\{11,12,13,$ $14,15,16,$ $17,18,19,$  $20,21,22,$ $23,58,59,$ $60,61,62,$ $63,64,65,$ $66,67,68,$ $69\}$, $X_4=\{24,25,26,$ $27,28,29,$ $30,30,31,$  $32,33,34,$ $35,$ $36,37,$ $38,39,40,$ $41,42,43,$ $45,70,71,$ $45,70,71,$ $72,73,74,$ $75,76,77,$ $78,79,80,$ $81,82,83,$  $84,85,86,$  $87,$ $88,89,$ $90,91\}$.}

\bigskip
\begin{center}

\hrule\vskip4pt

\end{center}

\bigskip

\footnotesize{
\noindent \textbf{BSTS($99$)}

\noindent  Color classes of Table~23: $X_1=\{1\}$, $X_2=\{2,3,4,$ $5\}$, $X_3=\{6,7,8,$ $9,50,$ $51,$ $52,53,$ $54,$  $55,56,57,$ $58,59,60,$ $61\}$, $X_4=\{10,11,12,$ $13,14,15,$ $16,17,18,$  $19,20,21,$ $22,23,24,$ $25,26,$ $27,$ $28,29,62,$ $63,64,65,$ $66,$ $67,68,$ $69,70,71,$ $72,$ $73,$ $74,$ $75,$ $76,77\}$ e $X_5=\{30,31,32,$ $33,34,35,$ $36,$ $37,38,$  $39,40,41,$ $42,43,44,$ $45,46,47,$ $48,49,$ $78,$ $79,80,81,$ $82,83,84,$ $85,86,$ $87,$ $88,$ $89,90,$ $91,$ $92,$ $93,$ $94,$ $95,96,$ $97,$ $98,$ $99\}$.}

\bigskip

\footnotesize{

\noindent  Color classes of Table~24: $X_1=\{1,2,$ $3$ $4,5,50,$ $51,$ $ 52\}$ $X_2=\{6,$ $7,8,$  $9,$ $ 10,$ $ 11$  $53,$ $ 54, $
$ 55$ $56,$ $ 57, $ $58$ $59, 60, 61$ $62,63\}$, $X_3=\{12,13,14,$ $15,16, 17$ $18, 19, 20$ $21, 22, 23$ $24, 25, 64$ $65, 66,$ $ 67$ $68, 69, 70$ $71,$ $ 72, 73$ $74, 75\}$, $X_4=\{26,$ $27,28,$ $29,30,31,$ $32,33,34,$  $35,36,37,$ $38,39,40,$ $41,42,43,$ $44,45,46,$ $47,48,49,$ $76,$ $77,78,$ $79,80,81,$ $82,83,84,$ $85,$ $86,87, 88,$ $89,90,$ $91,92,$ $ 93,$ $94,95,96,$  $97,98,$ $99\}$.}

\bigskip

\footnotesize{

\noindent  Color classes of Table~25: $X_1=\{1,2,50,$ $51, 52, $ $53, 54, 55\}$    $X_2=\{3,4,5$ $ 6,7,8,$  $9, 10,$ $56,57,$ $ 58,59 \}$, $X_3=\{11, 12,13,14,$ $15,16, 17$ $18, 19, 20$ $21, 22, 23$ $24, 25, $ $26,27,$ $28, 60,$ $61, 62, 63$ $64, 65, $ $66,$ $ 67$ $68, 69, 70$ $71, 72, 73$ $74, 75, 76$  $77, 78, 79\}$, $X_4=\{29,30,$ $31,$ $32,$ $33,34,$ $35,36,37,$  $38,39,40,$ $41,42,43,$ $44,45,46,$ $47,48,$ $49,$ $80,81,82,$ $83,$ $84,85,$ $86,87,88,$ $89,90,91,$ $92,$ $93,94, 95,$ $96,97,$ $98,99\}$.}

\bigskip

\footnotesize{
\noindent  Color classes of Table~26: $X_1=\{1,2,50,$ $51, 52, $ $53, 54, 55$ $56, $ $57, 58,59\}$ $X_2=\{3,$ $4,$ $5,$ $ 6,7,8,$  $9, 10\}$, $X_3=\{11, 12,13,14,$ $15,16, 17$ $18, 19, 20$ $21, 22, 23$ $24, 25, $ $26,27,$ $28, 60,$ $61, 62,$ $ 63$ $64, 65, $ $66, 67$ $68, 69, 70$ $71, 72, 73$ $74, 75, 76$  $77, 78, 79\}$, $X_4=\{29,30,$ $31,$ $32,33,$ $34,$ $35,36,37,$  $38,39,40,$ $41,42,43,$ $44,45,46,$ $47,48,$ $49,$ $80,81,82,$ $83,$ $84,85,$ $86,87,88,$ $89,90,91,$ $92,$ $93,94, 95,$ $96,97,$ $98,99\}$.}

\footnotesize{

\begin{center}

\begin{sideways}
\begin{tabular}{ccccccccccccccc}
~~&~~&~~&~~&~~&~~&~BSTS($91$)~&~~&~~&~~&~~&~~&~~&~~&~~\\
1&2&3&4&5&6&7&8&9&10&11&12&13&14&15\\
46~~50&46~~54&46~~51&46~~55&50~~83&50~~48&50~~72&50~~78&50~~86&50~~49&58~~46&58~~47&58~~48&58~~49&58~~57\\
47~~51&47~~55&47~~50&47~~52&51~~84&51~~49&51~~73&51~~70&51~~89&51~~48&59~~47&59~~46&59~~49&59~~48&59~~56\\
48~~52&48~~56&48~~53&48~~57&52~~85&52~~71&52~~74&52~~76&52~~46&52~~72&60~~48&60~~49&60~~46&60~~47&60~~50\\
49~~53&49~~57&49~~52&49~~56&53~~86&53~~81&53~~75&53~~77&53~~47&53~~46&61~~49&61~~48&61~~47&61~~46&61~~51\\
58~~59&58~~60&58~~61&58~~62&54~~87&54~~82&54~~79&54~~49&54~~48&54~~47&62~~50&62~~51&62~~52&62~~53&62~~46\\
60~~69&59~~61&60~~62&61~~63&55~~88&55~~89&55~~48&55~~87&55~~49&55~~71&63~~51&63~~50&63~~53&63~~52&63~~47\\
61~~68&62~~69&59~~63&60~~64&56~~46&56~~90&56~~80&56~~47&56~~91&56~~73&64~~52&64~~53&71~~72&64~~55&64~~48\\
62~~67&63~~68&64~~69&59~~65&57~~47&57~~91&57~~46&57~~88&57~~90&57~~85&65~~53&65~~52&64~~54&65~~54&65~~49\\
63~~66&64~~67&65~~68&66~~69&58~~63&58~~64&58~~65&58~~66&58~~67&58~~68&66~~54&66~~55&65~~55&66~~57&66~~52\\
64~~65&65~~66&66~~67&67~~68&62~~64&63~~65&64~~66&65~~67&66~~68&67~~69&67~~55&67~~54&66~~56&67~~56&67~~53\\
70~~71&70~~72&70~~73&70~~74&61~~65&62~~66&63~~67&64~~68&65~~69&59~~66&68~~56&68~~57&67~~57&68~~51&68~~54\\
72~~91&71~~73&72~~74&73~~75&60~~66&61~~67&62~~68&63~~69&59~~64&60~~65&69~~57&69~~56&68~~50&69~~50&69~~55\\
73~~90&74~~91&71~~75&72~~76&59~~67&60~~68&61~~69&59~~62&60~~63&61~~64&70~~80&70~~81&69~~51&70~~83&70~~84\\
74~~89&75~~90&76~~91&71~~77&68~~69&59~~69&59~~60&60~~61&61~~62&62~~63&79~~81&80~~82&70~~82&82~~84&83~~85\\
75~~88&76~~89&77~~90&78~~91&70~~75&83~~88&71~~81&75~~79&87~~88&86~~91&78~~82&79~~83&81~~83&81~~85&82~~86\\
76~~87&77~~88&78~~89&79~~90&74~~76&84~~87&82~~91&74~~80&70~~78&87~~90&77~~83&78~~84&80~~84&80~~86&81~~87\\
77~~86&78~~87&79~~88&80~~89&73~~77&85~~86&83~~90&73~~81&77~~79&88~~89&76~~84&77~~85&79~~85&79~~87&80~~88\\
78~~85&79~~86&80~~87&81~~88&72~~78&70~~76&84~~89&72~~82&76~~80&70~~79&75~~85&76~~86&78~~86&78~~88&79~~89\\
79~~84&80~~85&81~~86&82~~87&71~~79&75~~77&85~~88&71~~83&75~~81&78~~80&74~~86&75~~87&77~~87&77~~89&78~~90\\
80~~83&81~~84&82~~85&83~~86&80~~91&74~~78&86~~87&84~~91&74~~82&77~~81&73~~87&74~~88&76~~88&76~~90&77~~91\\
81~~82&82~~83&83~~84&84~~85&81~~90&73~~79&70~~77&85~~90&73~~83&76~~82&72~~88&73~~89&75~~89&75~~91&71~~76\\
54~~57&50~~52&54~~56&50~~51&82~~89&72~~80&76~~78&86~~89&72~~84&75~~83&71~~89&72~~90&74~~90&71~~74&72~~75\\
55~~56&51~~53&55~~57&53~~54&48~~49&46~~47&47~~49&46~~48&71~~85&74~~84&90~~91&71~~91&73~~91&72~~73&73~~74\\
~~&~~&~~&~~&~~&~~&~~&~~&~~&~~&~~&~~&~~&~~&~~\\
~~&~~&~~&~~&~~&~~&table~22&~~&~&~&~&~&~&~~&\\
\end{tabular}

\end{sideways}

\end{center}}

\footnotesize{
\begin{center}

\begin{sideways}
\begin{tabular}{ccccccccccccccc}
~~&~~&~~&~~&~~&~~&~BSTS($91$)~&~~&~~&~~&~~&~~&~~&~~&~~\\
16&17&18&19&20&21&22&23&24&25&26&27&28&29&30\\
58~~56&58~~54&58~~51&58~~53&58~~52&58~~55&58~~50&58~~78&70~~59&70~~61&70~~58&70~~55&70~~60&70~~62&70~~50\\
59~~57&59~~55&59~~52&59~~50&59~~53&59~~54&59~~51&59~~80&71~~60&71~~59&71~~61&71~~62&71~~58&71~~63&71~~51\\
60~~51&60~~56&60~~53&60~~52&60~~54&60~~57&60~~74&60~~55&72~~61&72~~58&72~~57&72~~60&72~~62&72~~64&72~~54\\
61~~50&61~~57&61~~55&61~~54&61~~56&61~~52&61~~53&61~~82&73~~62&73~~63&73~~67&73~~59&73~~61&73~~65&73~~55\\
62~~47&62~~48&62~~54&62~~55&62~~57&62~~49&62~~56&62~~79&74~~63&74~~62&74~~66&74~~54&74~~64&74~~67&74~~56\\
63~~46&63~~49&63~~57&63~~56&63~~55&63~~48&63~~70&63~~54&75~~64&75~~66&75~~60&75~~61&75~~63&75~~68&75~~57\\
64~~49&64~~46&64~~56&64~~57&64~~50&64~~51&64~~47&64~~81&76~~65&76~~64&76~~68&76~~56&76~~67&76~~69&76~~47\\
65~~48&65~~47&65~~50&65~~51&65~~46&65~~56&65~~75&65~~57&77~~66&77~~69&77~~63&77~~67&77~~65&77~~46&77~~48\\
66~~53&66~~50&66~~49&66~~48&66~~51&66~~47&66~~76&66~~46&78~~67&78~~65&78~~64&78~~48&78~~66&78~~47&78~~46\\
67~~52&67~~51&67~~46&67~~49&67~~47&67~~50&67~~48&67~~83&79~~68&79~~46&79~~59&79~~66&79~~69&79~~48&79~~49\\
68~~55&68~~52&68~~47&68~~46&68~~48&68~~53&68~~77&68~~49&80~~46&80~~68&80~~51&80~~69&80~~47&80~~49&80~~67\\
69~~54&69~~53&69~~48&69~~47&69~~49&69~~46&69~~84&69~~52&81~~47&81~~48&81~~54&81~~68&81~~46&81~~55&81~~69\\
70~~85&70~~86&70~~87&70~~88&70~~89&70~~90&73~~85&70~~91&82~~48&82~~47&82~~46&82~~65&82~~49&82~~56&82~~68\\
84~~86&85~~87&86~~88&87~~89&88~~90&89~~91&72~~86&71~~90&83~~49&83~~51&83~~55&83~~53&83~~48&83~~57&83~~66\\
83~~87&84~~88&85~~89&86~~90&87~~91&71~~88&71~~87&72~~89&84~~50&84~~49&84~~53&84~~57&84~~52&84~~58&84~~65\\
82~~88&83~~89&84~~90&85~~91&71~~86&72~~87&88~~91&73~~88&85~~51&85~~56&85~~69&85~~58&85~~50&85~~53&85~~64\\
81~~89&82~~90&83~~91&71~~84&72~~85&73~~86&89~~90&74~~87&86~~52&86~~55&86~~56&86~~51&86~~57&86~~59&86~~63\\
80~~90&81~~91&71~~82&72~~83&73~~84&74~~85&80~~81&75~~86&87~~57&87~~50&87~~52&87~~47&87~~53&87~~60&87~~62\\
79~~91&71~~80&72~~81&73~~82&74~~83&75~~84&79~~82&76~~85&88~~54&88~~52&88~~48&88~~49&88~~56&88~~51&88~~61\\
71~~78&72~~79&73~~80&74~~81&75~~82&76~~83&78~~83&77~~84&89~~56&89~~57&89~~49&89~~52&89~~54&89~~50&89~~60\\
72~~77&73~~78&74~~79&75~~80&76~~81&77~~82&46~~49&47~~48&90~~53&90~~54&90~~50&90~~46&90~~55&90~~52&90~~59\\
73~~76&74~~77&75~~78&76~~79&77~~80&78~~81&52~~57&50~~53&91~~55&91~~53&91~~47&91~~50&91~~51&91~~54&91~~58\\
74~~75&75~~76&76~~77&77~~78&78~~79&79~~80&54~~55&51~~56&58~~69&60~~67&62~~65&63~~64&59~~68&61~~66&52~~53\\
~~&~~&~~&~~&~~&~~&~~&~~&~~&~~&~~&~~&~~&~~&~~\\
~~&~~&~~&~~&~~&~~&table~22&~~&~&~&~&~&~&~~&\\
\end{tabular}

\end{sideways}

\end{center}}

\footnotesize{
\begin{center}
\begin{sideways}
\begin{tabular}{ccccccccccccccc}
31&32&33&34&35&36&37&38&39&40&41&42&43&44&45\\
70~~52&70~~46&70~~54&70~~68&70~~47&70~~56&70~~69&70~~57&70~~66&70~~53&70~~64&70~~48&70~~49&70~~65&70~~67\\
71~~53&71~~47&71~~50&71~~69&71~~46&71~~68&71~~65&71~~56&71~~67&71~~54&71~~49&71~~66&71~~48&71~~57&71~~64\\
72~~55&72~~48&72~~51&72~~67&72~~49&72~~53&72~~66&72~~59&72~~68&72~~56&72~~46&72~~47&72~~69&72~~63&72~~65\\
73~~54&73~~49&73~~52&73~~60&73~~48&73~~58&73~~64&73~~53&73~~50&73~~69&73~~68&73~~57&73~~46&73~~66&73~~47\\
74~~50&74~~51&74~~57&74~~61&74~~53&74~~48&74~~68&74~~69&74~~59&74~~46&74~~65&74~~55&74~~47&74~~49&74~~58\\
75~~46&75~~52&75~~56&75~~50&75~~54&75~~47&75~~67&75~~49&75~~51&75~~48&75~~55&75~~58&75~~59&75~~62&75~~69\\
76~~48&76~~53&76~~58&76~~51&76~~50&76~~54&76~~63&76~~60&76~~61&76~~55&76~~57&76~~49&76~~62&76~~59&76~~46\\
77~~47&77~~55&77~~59&77~~57&77~~52&77~~61&77~~62&77~~64&77~~58&77~~50&77~~51&77~~56&77~~60&77~~54&77~~49\\
78~~49&78~~56&78~~62&78~~53&78~~55&78~~63&78~~61&78~~68&78~~52&78~~57&78~~54&78~~59&78~~51&78~~69&78~~60\\
79~~51&79~~57&79~~60&79~~55&79~~56&79~~64&79~~58&79~~65&79~~63&79~~47&79~~52&79~~50&79~~53&79~~67&79~~61\\
80~~66&80~~60&80~~61&80~~54&80~~58&80~~65&80~~57&80~~48&80~~62&80~~64&80~~53&80~~52&80~~55&80~~50&80~~63\\
81~~67&81~~61&81~~63&81~~62&81~~59&81~~57&81~~52&81~~51&81~~56&81~~49&81~~60&81~~65&81~~66&81~~58&81~~50\\
82~~69&82~~66&82~~67&82~~58&82~~60&82~~59&82~~51&82~~50&82~~55&82~~62&82~~63&82~~53&82~~52&82~~64&82~~57\\
83~~68&83~~63&83~~64&83~~65&83~~61&83~~46&83~~60&83~~47&83~~54&83~~58&83~~62&83~~69&83~~56&83~~52&83~~59\\
84~~64&84~~67&84~~66&84~~63&84~~62&84~~60&84~~59&84~~55&84~~48&84~~68&84~~61&84~~46&84~~54&84~~47&84~~56\\
85~~65&85~~62&85~~68&85~~66&85~~63&85~~67&85~~55&85~~46&85~~49&85~~59&85~~47&85~~60&85~~61&85~~48&85~~54\\
86~~62&86~~65&86~~69&86~~64&86~~67&86~~49&86~~54&86~~66&86~~47&86~~60&86~~48&86~~61&86~~58&86~~46&86~~68\\
87~~63&87~~64&87~~65&87~~59&87~~69&87~~51&87~~49&87~~58&87~~46&87~~61&87~~66&87~~68&87~~67&87~~56&87~~48\\
88~~60&88~~58&88~~47&88~~46&88~~64&88~~66&88~~50&88~~67&88~~69&88~~65&88~~59&88~~62&88~~63&88~~68&88~~53\\
89~~61&89~~59&89~~46&89~~47&89~~65&89~~69&89~~48&89~~62&89~~64&89~~63&89~~58&89~~67&89~~68&89~~53&89~~66\\
90~~58&90~~69&90~~48&90~~49&90~~68&90~~62&90~~47&90~~61&90~~65&90~~66&90~~67&90~~63&90~~64&90~~60&90~~51\\
91~~59&91~~68&91~~49&91~~48&91~~66&91~~52&91~~46&91~~63&91~~60&91~~67&91~~69&91~~64&91~~65&91~~61&91~~62\\
56~~57&50~~54&53~~55&52~~56&51~~57&50~~55&53~~56&52~~54&53~~57&51~~52&50~~56&51~~54&50~~57&51~~55&52~~55\\
~~&~~&~~&~~&~~&~~&~~&~~&~~&~~&~~&~~&~~&~~&~~\\
~~&~~&~~&~~&~~&~~&~~&~table~22~&~&~&~&~&~&~~&\\
\end{tabular}

\end{sideways}

\end{center}}

\footnotesize{

\begin{center}
\begin{sideways}
\begin{tabular}{ccccccccccccccccc}  \label{tbl1}
~~&~~&~~&~~&~~&~~&~~&BSTS($99$)&~~&~~&~~&~~&~~&~~&~~&~~&~~\\
1&2&3&4&5&6&7&8&9&10&11&12&13&14&15&16&17\\
62~~67&62~~68&62~~72&62~~73&62~~75&50~~62&50~~69&50~~88&50~~90&62~~51&62~~52&62~~53&62~~54&62~~55&62~~99&62~~57&62~~58\\
63~~71&63~~73&63~~75&63~~65&63~~72&51~~76&51~~85&51~~92&51~~66&63~~52&64~~53&63~~55&63~~53&63~~91&63~~57&63~~58&63~~78\\
64~~77&64~~72&64~~68&64~~69&64~~65&52~~70&52~~77&52~~86&52~~88&64~~50&65~~61&64~~51&64~~52&64~~54&64~~58&64~~80&64~~57\\
65~~74&65~~67&65~~76&66~~68&66~~69&53~~74&53~~93&53~~97&53~~73&65~~60&66~~50&65~~52&65~~58&65~~57&65~~59&65~~84&65~~50\\
66~~75&66~~71&66~~74&67~~72&67~~77&54~~77&54~~95&54~~63&54~~98&66~~53&67~~55&66~~61&66~~60&66~~59&66~~52&66~~55&66~~54\\
68~~70&70~~76&67~~71&71~~77&68~~74&55~~97&55~~75&55~~96&55~~64&67~~54&68~~60&67~~60&67~~61&67~~94&67~~53&67~~52&67~~56\\
69~~73&69~~77&69~~70&75~~76&73~~76&56~~95&56~~97&56~~68&56~~74&68~~61&69~~59&68~~59&68~~50&68~~51&68~~54&68~~53&68~~52\\
72~~76&74~~75&73~~77&70~~74&70~~71&57~~78&57~~74&57~~71&57~~93&69~~55&70~~57&69~~58&69~~83&69~~60&69~~61&69~~54&69~~51\\
78~~84&78~~86&78~~87&78~~88&78~~91&58~~72&58~~78&58~~74&58~~92&70~~58&71~~89&70~~56&70~~55&70~~95&70~~88&70~~51&70~~91\\
79~~98&79~~80&79~~94&79~~90&79~~99&59~~92&59~~73&59~~93&59~~76&71~~56&72~~54&71~~54&71~~59&71~~50&71~~55&71~~60&71~~61\\
80~~91&81~~95&80~~98&80~~95&80~~88&60~~88&60~~90&60~~76&60~~77&72~~59&73~~56&72~~57&72~~56&72~~96&72~~50&72~~61&72~~60\\
81~~93&82~~92&81~~91&81~~89&85~~93&61~~84&61~~76&61~~77&61~~87&73~~57&74~~80&73~~50&73~~51&73~~52&73~~97&73~~93&73~~90\\
82~~96&83~~90&82~~97&82~~86&82~~94&80~~99&79~~81&78~~79&78~~99&74~~91&75~~58&74~~81&74~~90&74~~61&74~~51&74~~50&74~~84\\
83~~85&84~~94&83~~89&83~~94&83~~95&81~~98&80~~84&80~~82&79~~83&75~~87&76~~86&75~~89&75~~57&75~~53&75~~60&75~~56&75~~59\\
86~~90&85~~96&84~~93&87~~93&86~~98&83~~96&82~~99&81~~83&80~~97&76~~81&77~~92&76~~92&76~~99&76~~58&76~~56&76~~94&76~~53\\
87~~97&87~~99&85~~99&84~~92&89~~96&85~~94&83~~98&84~~95&81~~96&77~~82&78~~81&77~~78&77~~98&77~~56&77~~79&77~~59&77~~55\\
88~~94&88~~91&90~~92&85~~98&84~~90&86~~93&86~~88&85~~89&82~~95&78~~80&79~~87&84~~99&78~~97&78~~98&78~~94&78~~96&79~~89\\
89~~92&89~~93&86~~96&91~~99&87~~92&89~~90&87~~89&87~~94&84~~91&79~~85&82~~93&83~~87&79~~84&93~~99&80~~81&79~~86&80~~83\\
95~~99&97~~98&88~~95&96~~97&81~~97&79~~82&91~~94&90~~91&85~~86&83~~88&83~~99&82~~88&80~~85&79~~97&82~~83&81~~85&81~~92\\
50~~51&50~~52&50~~53&50~~54&50~~55&87~~91&92~~96&98~~99&89~~94&84~~86&84~~98&80~~90&81~~82&80~~86&84~~89&82~~90&82~~85\\
52~~61&51~~53&52~~54&53~~55&54~~56&63~~64&62~~65&62~~70&62~~71&89~~99&85~~88&79~~91&86~~87&82~~89&85~~91&83~~91&94~~97\\
53~~60&54~~61&51~~55&52~~56&53~~57&65~~73&63~~67&64~~67&63~~70&90~~98&90~~96&93~~97&88~~89&81~~84&86~~95&87~~95&86~~99\\
54~~59&55~~60&56~~61&51~~57&52~~58&66~~67&64~~71&65~~69&65~~75&92~~94&94~~95&86~~94&91~~95&87~~88&87~~96&88~~98&87~~98\\
55~~58&56~~59&57~~60&58~~61&51~~59&68~~75&66~~70&66~~73&67~~68&93~~96&91~~97&85~~95&92~~93&85~~90&90~~93&89~~97&88~~93\\
56~~57&57~~58&58~~59&59~~60&60~~61&69~~71&68~~72&72~~75&69~~72&95~~97&63~~51&96~~98&94~~96&83~~92&92~~98&92~~99&95~~96\\
~~&~~&~~&~~&~~&~~&~~&~~&~~&~~&~~&~~&~~&~~&~~&~~&~~\\
~~&~~&~~&~~&~~&~~&~~&table~23&~~&~&~&~&~&~&~~&\\
\end{tabular}

\end{sideways}
\end{center}}

\footnotesize{
\begin{center}

\begin{sideways}
\begin{tabular}{cccccccccccccccc}
18&19&20&21&22&23&24&25&26&27&28&29&30&31&32&33\\
62~~98&62~~60&62~~61&62~~59&62~~79&62~~56&62~~80&62~~82&62~~84&62~~89&62~~91&62~~94&78~~60&78~~62&78~~64&78~~65\\
63~~60&63~~61&63~~50&63~~56&63~~59&63~~82&63~~79&63~~83&63~~86&63~~88&63~~97&63~~87&79~~64&79~~76&79~~66&79~~71\\
64~~96&64~~56&64~~99&64~~61&64~~60&64~~59&64~~95&64~~92&64~~90&64~~83&64~~81&64~~85&80~~70&80~~65&80~~63&80~~72\\
65~~54&65~~55&65~~53&65~~79&65~~56&65~~51&65~~81&65~~85&65~~91&65~~94&65~~98&65~~99&81~~62&81~~63&81~~67&81~~66\\
66~~56&66~~85&66~~58&66~~80&66~~57&66~~87&66~~83&66~~84&66~~98&66~~96&66~~92&66~~95&82~~75&82~~64&82~~65&82~~60\\
67~~57&67~~50&67~~59&67~~58&67~~51&67~~89&67~~99&67~~98&67~~95&67~~80&67~~79&67~~84&83~~50&83~~74&83~~52&83~~53\\
68~~55&68~~93&68~~96&68~~86&68~~58&68~~57&68~~94&68~~88&68~~79&68~~82&68~~83&68~~80&84~~54&84~~50&84~~51&84~~52\\
69~~52&69~~53&69~~56&69~~57&69~~85&69~~91&69~~82&69~~80&69~~81&69~~97&69~~90&69~~88&85~~53&85~~52&85~~60&85~~59\\
70~~50&70~~59&70~~54&70~~53&70~~61&70~~60&70~~84&70~~81&70~~93&70~~99&70~~82&70~~83&86~~56&86~~59&86~~50&86~~51\\
71~~53&71~~51&71~~52&71~~93&71~~86&71~~58&71~~97&71~~91&71~~87&71~~84&71~~85&71~~81&87~~59&87~~67&87~~73&87~~64\\
72~~51&72~~52&72~~94&72~~55&72~~89&72~~53&72~~87&72~~79&72~~99&72~~85&72~~78&72~~82&88~~55&88~~51&88~~53&88~~61\\
73~~58&73~~54&73~~95&73~~60&73~~55&73~~61&73~~98&73~~96&73~~82&73~~81&73~~80&73~~79&89~~61&89~~66&89~~70&89~~73\\
74~~59&74~~98&74~~60&74~~52&74~~54&74~~55&74~~88&74~~87&74~~89&74~~78&74~~86&74~~96&90~~58&90~~53&90~~56&90~~54\\
75~~61&75~~95&75~~51&75~~50&75~~52&75~~54&75~~90&75~~99&75~~78&75~~86&75~~94&75~~92&91~~73&91~~55&91~~57&91~~56\\
76~~93&76~~57&76~~55&76~~54&76~~50&76~~52&76~~89&76~~78&76~~83&76~~87&76~~96&76~~91&92~~67&92~~60&92~~54&92~~55\\
77~~83&77~~58&77~~57&77~~51&77~~53&77~~50&77~~85&77~~86&77~~94&77~~90&77~~89&77~~97&93~~77&93~~61&93~~69&93~~62\\
78~~89&78~~92&78~~83&78~~95&78~~82&78~~85&78~~93&89~~95&80~~96&79~~92&84~~87&78~~90&94~~57&94~~56&94~~55&94~~58\\
79~~95&79~~96&79~~93&85~~87&81~~87&79~~88&86~~92&90~~97&88~~97&91~~93&88~~99&86~~89&95~~65&95~~69&95~~59&95~~68\\
80~~94&80~~89&80~~87&90~~94&80~~92&80~~93&91~~96&93~~94&85~~92&95~~98&93~~95&93~~98&96~~52&96~~54&96~~58&96~~57\\
81~~88&81~~94&81~~86&89~~91&83~~93&81~~90&50~~56&50~~57&50~~58&50~~59&50~~60&50~~61&97~~51&97~~58&97~~61&97~~50\\
82~~87&82~~91&82~~84&84~~96&84~~97&83~~84&55~~57&56~~58&57~~59&58~~60&59~~61&51~~60&98~~71&98~~68&98~~72&98~~63\\
84~~85&83~~86&85~~97&83~~97&88~~96&86~~97&54~~58&55~~59&56~~60&57~~61&51~~58&52~~59&99~~68&99~~57&99~~77&99~~74\\
86~~91&84~~88&88~~90&82~~98&90~~95&92~~95&53~~59&54~~60&55~~61&51~~56&52~~57&53~~58&63~~66&70~~72&68~~76&67~~70\\
90~~99&87~~90&91~~92&81~~99&94~~99&94~~98&52~~60&53~~61&51~~54&52~~55&53~~56&54~~57&69~~76&71~~73&71~~75&69~~75\\
92~~97&97~~99&89~~98&88~~92&91~~98&96~~99&51~~61&51~~52&52~~53&53~~54&54~~55&55~~56&72~~74&75~~77&62~~74&76~~77\\
~~&~~&~~&~~&~~&~~&~~&~~&~~&~~&~~&~~&~~&~~&~~&~~\\
~~&~~&~~&~~&~~&~~&~~&~~&~table~23~&~&~&~&~&~&~~&\\
\end{tabular}

\end{sideways}
\end{center}}

\footnotesize{
\begin{center}
\begin{sideways}
\begin{tabular}{cccccccccccccccc}
34&35&36&37&38&39&40&41&42&43&44&45&46&47&48&49\\
78~~66&78~~59&78~~68&78~~69&78~~70&78~~71&78~~73&78~~50&78~~51&78~~52&78~~61&78~~54&78~~55&78~~67&78~~56&78~~53\\
79~~75&79~~51&79~~52&79~~53&79~~54&79~~55&79~~56&79~~57&79~~58&79~~59&79~~60&79~~61&79~~74&79~~69&79~~50&79~~70\\
80~~71&80~~76&80~~77&80~~50&80~~51&80~~52&80~~53&80~~54&80~~55&80~~56&80~~57&80~~58&80~~59&80~~60&80~~61&80~~75\\
81~~68&81~~72&81~~75&81~~61&81~~50&81~~51&81~~52&81~~77&81~~54&81~~55&81~~53&81~~57&81~~58&81~~59&81~~60&81~~56\\
82~~51&82~~66&82~~67&82~~71&82~~52&82~~53&82~~54&82~~55&82~~61&82~~76&82~~56&82~~50&82~~57&82~~58&82~~74&82~~59\\
83~~54&83~~55&83~~56&83~~57&83~~58&83~~62&83~~59&83~~60&83~~72&83~~61&83~~73&83~~75&83~~67&83~~51&83~~65&83~~71\\
84~~53&84~~60&84~~55&84~~59&84~~75&84~~58&84~~57&84~~56&84~~63&84~~77&84~~69&84~~76&84~~68&84~~64&84~~73&84~~72\\
85~~58&85~~57&85~~54&85~~56&85~~55&85~~63&85~~50&85~~61&85~~76&85~~67&85~~68&85~~62&85~~70&85~~73&85~~75&85~~74\\
86~~61&86~~69&86~~53&86~~66&86~~67&86~~73&86~~70&86~~72&86~~60&86~~62&86~~54&86~~55&86~~65&86~~57&86~~58&86~~64\\
87~~65&87~~68&87~~70&87~~58&87~~53&87~~50&87~~77&87~~52&87~~57&87~~54&87~~62&87~~69&87~~56&87~~55&87~~51&87~~60\\
88~~59&88~~58&88~~57&88~~76&88~~56&88~~77&88~~75&88~~62&88~~64&88~~65&88~~67&88~~66&88~~71&88~~54&88~~72&88~~73\\
89~~64&89~~65&89~~63&89~~68&89~~57&89~~69&89~~58&89~~59&89~~56&89~~60&89~~50&89~~51&89~~53&89~~52&89~~55&89~~54\\
90~~55&90~~52&90~~51&90~~63&90~~62&90~~68&90~~66&90~~65&90~~67&90~~57&90~~70&90~~71&90~~72&90~~61&90~~59&90~~76\\
91~~52&91~~50&91~~60&91~~54&91~~77&91~~59&91~~64&91~~51&91~~66&91~~53&91~~75&91~~68&91~~61&91~~72&91~~67&91~~58\\
92~~56&92~~61&92~~50&92~~62&92~~63&92~~57&92~~65&92~~69&92~~52&92~~70&92~~71&92~~72&92~~73&92~~74&92~~53&92~~68\\
93~~63&93~~56&93~~64&93~~65&93~~66&93~~67&93~~72&93~~74&93~~75&93~~58&93~~51&93~~60&93~~52&93~~50&93~~54&93~~55\\
94~~60&94~~54&94~~59&94~~52&94~~73&94~~74&94~~61&94~~71&94~~53&94~~51&94~~66&94~~70&94~~64&94~~63&94~~69&94~~50\\
95~~62&95~~63&95~~71&95~~72&95~~74&95~~76&95~~51&95~~58&95~~77&95~~50&95~~55&95~~52&95~~60&95~~53&95~~57&95~~61\\
96~~50&96~~53&96~~61&96~~51&96~~60&96~~56&96~~62&96~~63&96~~65&96~~69&96~~59&96~~67&96~~75&96~~71&96~~70&96~~77\\
97~~57&97~~62&97~~65&97~~64&97~~72&97~~60&97~~67&97~~68&97~~70&97~~74&97~~76&97~~59&97~~54&97~~75&97~~66&97~~52\\
98~~76&98~~64&98~~69&98~~55&98~~59&98~~61&98~~60&98~~53&98~~50&98~~75&98~~58&98~~56&98~~51&98~~70&98~~52&98~~57\\
99~~69&99~~71&99~~58&99~~60&99~~61&99~~54&99~~55&99~~66&99~~59&99~~73&99~~52&99~~53&99~~50&99~~56&99~~63&99~~51\\
67~~74&67~~73&62~~66&67~~75&64~~76&64~~75&63~~76&64~~70&62~~69&63~~68&63~~77&63~~74&62~~76&62~~77&62~~64&62~~63\\
70~~73&70~~75&72~~73&70~~77&65~~71&65~~70&68~~71&67~~76&68~~73&64~~66&64~~74&65~~77&63~~69&65~~68&68~~77&65~~66\\
72~~77&74~~77&74~~76&73~~74&68~~69&66~~72&69~~74&73~~75&71~~74&71~~72&65~~72&64~~73&66~~77&66~~76&71~~76&67~~69\\
~~&~~&~~&~~&~~&~~&~~&~~&~~&~~&~~&~~&~~&~~&~~&~~\\
~~&~~&~~&~~&~~&~~&~~&~~&~table~23~&~&~&~&~&~&~~&\\
\end{tabular}

\end{sideways}
\end{center}}

\footnotesize{
\begin{center}
\begin{sideways}

\begin{tabular}{ccccccccccccccccc}  \label{tbl1}
~~&~~&~~&~~&~~&~~&~~&~\textbf{BSTS($99$)}~&~~&~~&~~&~~&~~&~~&~~&~~&~~\\
1&2&3&4&5&6&7&8&9&10&11&12&13&14&15&16&17\\
50~~54&50~~55&50~~56&50~~57&50~~58&59~~50&60~~50&61~~50&62~~50&63~~50&53~~50&64~~50&65~~50&66~~50&67~~50&67~~52&69~~50\\
51~~53&51~~54&51~~55&51~~56&51~~57&58~~51&59~~51&60~~51&61~~51&62~~51&63~~51&70~~52&74~~51&73~~51&72~~51&68~~50&70~~51\\
52~~56&52~~57&52~~58&52~~59&52~~60&61~~52&62~~52&63~~52&53~~52&54~~52&55~~52&75~~51&69~~52&68~~52&66~~52&71~~51&65~~52\\
55~~57&53~~58&54~~53&63~~55&54~~59&53~~68&53~~64&53~~66&54~~64&53~~65&54~~75&65~~80&64~~55&64~~60&64~~63&64~~62&64~~58\\
58~~59&56~~59&60~~57&61~~53&53~~63&54~~72&54~~71&54~~67&55~~73&55~~72&56~~64&66~~96&66~~59&65~~59&65~~77&65~~60&66~~61\\
60~~61&60~~63&61~~63&62~~60&62~~56&55~~67&55~~68&55~~70&56~~65&56~~66&57~~65&67~~62&67~~79&67~~53&68~~56&66~~54&67~~63\\
62~~63&61~~62&59~~62&54~~58&61~~55&56~~73&56~~74&56~~71&57~~70&57~~71&58~~74&68~~61&68~~58&69~~54&69~~57&69~~63&68~~82\\
64~~65&64~~66&64~~67&64~~68&64~~69&57~~66&57~~67&57~~74&58~~66&58~~67&59~~67&69~~58&70~~62&70~~63&70~~84&70~~61&71~~83\\
66~~75&67~~65&66~~68&67~~69&68~~70&60~~74&58~~75&58~~73&59~~69&59~~70&60~~66&71~~55&71~~60&71~~62&71~~58&72~~94&72~~56\\
67~~74&68~~75&65~~69&66~~70&67~~71&62~~65&61~~65&59~~64&60~~67&60~~68&61~~72&72~~59&72~~57&72~~81&73~~60&73~~92&73~~57\\
68~~73&69~~74&70~~75&65~~71&66~~72&63~~75&63~~66&62~~72&63~~68&61~~69&62~~73&73~~54&73~~61&74~~78&74~~61&74~~55&74~~54\\
69~~72&70~~73&71~~74&72~~75&65~~73&64~~70&70~~72&69~~75&72~~74&64~~74&68~~71&74~~53&75~~99&75~~56&75~~59&75~~53&75~~62\\
71~~70&72~~71&72~~73&74~~73&75~~74&69~~71&69~~73&65~~68&71~~75&73~~75&69~~70&56~~60&53~~56&55~~58&53~~62&56~~58&53~~59\\
76~~77&76~~78&76~~79&76~~80&76~~81&76~~82&76~~83&76~~84&76~~85&76~~86&76~~87&57~~63&54~~63&57~~61&54~~55&59~~57&55~~60\\
78~~99&79~~77&78~~80&79~~81&80~~82&81~~83&82~~84&83~~85&84~~86&85~~87&86~~88&76~~88&76~~89&76~~91&76~~92&76~~93&76~~94\\
79~~98&80~~99&77~~81&78~~82&79~~83&80~~84&81~~85&82~~86&83~~87&84~~88&85~~89&87~~89&88~~90&90~~92&91~~93&91~~95&93~~95\\
80~~97&81~~98&82~~99&77~~83&78~~84&79~~85&80~~86&81~~87&82~~88&83~~89&84~~90&86~~90&87~~91&89~~93&90~~94&90~~96&92~~96\\
81~~96&82~~97&83~~98&84~~99&77~~85&78~~86&79~~87&80~~88&81~~89&82~~90&83~~91&85~~91&86~~92&88~~94&89~~95&89~~97&91~~97\\
82~~95&83~~96&84~~97&85~~98&86~~99&87~~77&88~~78&89~~79&90~~80&81~~91&82~~92&84~~92&85~~93&87~~95&88~~96&88~~98&90~~98\\
83~~94&84~~95&85~~96&86~~97&87~~98&88~~99&77~~89&78~~90&79~~91&80~~92&81~~93&93~~83&94~~84&86~~96&87~~97&87~~99&89~~99\\
84~~93&85~~94&86~~95&87~~96&88~~97&89~~98&90~~99&77~~91&78~~92&79~~93&80~~94&82~~94&83~~95&85~~97&86~~98&77~~86&77~~88\\
85~~92&86~~93&87~~94&88~~95&89~~96&90~~97&91~~98&92~~99&77~~93&78~~94&79~~95&81~~95&82~~96&84~~98&85~~99&85~~78&78~~87\\
86~~91&87~~92&88~~93&89~~94&90~~95&91~~96&92~~97&93~~98&94~~99&77~~95&78~~96&79~~97&81~~97&83~~99&78~~83&79~~84&79~~86\\
87~~90&88~~91&89~~92&90~~93&91~~94&92~~95&93~~96&94~~97&95~~98&96~~99&77~~97&78~~98&80~~98&77~~82&79~~82&80~~83&80~~85\\
88~~89&89~~90&90~~91&91~~92&92~~93&93~~94&94~~95&95~~96&96~~97&97~~98&98~~99&77~~99&77~~78&79~~80&80~~81&81~~82&81~~84\\
~~&~~&~~&~~&~~&~~&~~&~~&~~&~~&~~&~~&~~&~~&~~&~~&~~\\
~~&~~&~~&~~&~~&~~&~~&table~24&~~&~&~&~&~&~&~~&\\
\end{tabular}

\end{sideways}
\end{center}}

\footnotesize{
\begin{center}
\begin{sideways}

\begin{tabular}{cccccccccccccccc}
18&19&20&21&22&23&24&25&26&27&28&29&30&31&32&33\\
70~~50&71~~50&72~~50&73~~50&74~~50&50~~51&50~~52&75~~50&77~~50&76~~50&76~~52&77~~52&77~~56&76~~54&77~~61&76~~60\\
69~~51&68~~51&67~~51&66~~51&65~~51&71~~52&64~~51&64~~76&79~~52&79~~54&77~~54&78~~54&78~~52&78~~63&80~~54&78~~59\\
64~~52&75~~52&74~~52&72~~52&73~~52&64~~90&65~~90&65~~99&82~~54&80~~52&78~~50&79~~50&80~~50&80~~62&81~~52&79~~61\\
65~~54&64~~57&64~~85&64~~61&64~~91&65~~83&66~~85&66~~98&83~~56&82~~60&81~~60&81~~63&81~~54&81~~50&83~~50&82~~50\\
66~~80&65~~55&65~~63&65~~58&66~~62&66~~92&67~~86&67~~77&85~~58&84~~56&82~~62&84~~62&83~~63&83~~52&84~~57&84~~52\\
67~~56&66~~88&66~~55&67~~94&67~~89&67~~95&68~~93&68~~57&86~~63&86~~58&84~~63&87~~61&85~~62&85~~60&85~~56&86~~54\\
68~~87&67~~61&68~~62&68~~59&68~~54&68~~84&69~~95&69~~97&89~~61&87~~62&85~~61&88~~56&86~~61&87~~59&88~~63&87~~57\\
71~~61&69~~62&69~~60&69~~53&69~~55&69~~56&70~~76&70~~58&91~~62&89~~63&88~~58&89~~58&87~~60&88~~61&90~~62&89~~62\\
72~~53&70~~53&70~~54&70~~60&70~~56&70~~96&71~~88&71~~78&92~~60&91~~61&90~~59&90~~57&89~~59&89~~56&91~~60&90~~63\\
73~~63&72~~58&71~~59&71~~63&71~~53&72~~89&72~~97&72~~60&94~~59&93~~57&91~~57&91~~59&91~~58&92~~53&93~~55&91~~53\\
74~~62&73~~81&73~~53&74~~79&72~~63&73~~91&73~~87&73~~59&95~~57&95~~55&92~~55&93~~53&93~~51&93~~58&95~~51&92~~51\\
75~~60&74~~59&75~~86&75~~55&75~~57&74~~82&74~~98&74~~63&96~~55&96~~59&96~~56&94~~51&94~~53&96~~51&96~~53&94~~55\\
57~~58&56~~63&58~~61&54~~56&58~~60&75~~93&75~~61&56~~61&97~~53&98~~53&97~~51&97~~60&97~~57&97~~55&99~~58&96~~58\\
55~~59&54~~60&56~~57&57~~62&59~~61&53~~57&53~~60&53~~55&98~~51&99~~51&99~~53&98~~55&99~~55&99~~57&97~~59&97~~56\\
76~~95&76~~96&76~~97&76~~98&76~~90&54~~61&54~~57&54~~62&65~~76&65~~78&65~~93&65~~85&64~~98&64~~77&64~~78&64~~99\\
94~~96&95~~97&96~~98&97~~99&88~~92&55~~62&55~~56&79~~96&66~~87&66~~81&66~~95&66~~76&67~~76&66~~86&66~~89&66~~77\\
93~~97&94~~98&95~~99&96~~77&87~~93&58~~63&58~~62&80~~95&67~~99&67~~88&67~~98&67~~96&68~~88&67~~91&67~~82&67~~83\\
92~~98&93~~99&77~~94&78~~95&86~~94&59~~60&59~~63&81~~94&68~~81&68~~97&68~~86&68~~92&69~~84&68~~95&68~~98&68~~80\\
91~~99&77~~92&78~~93&80~~93&85~~95&76~~99&82~~83&82~~93&69~~90&69~~94&69~~79&70~~83&70~~79&69~~98&69~~76&69~~88\\
77~~90&78~~91&79~~92&81~~92&84~~96&77~~98&92~~94&83~~92&70~~93&70~~77&70~~87&71~~80&71~~82&71~~79&70~~86&70~~81\\
78~~89&79~~90&80~~91&82~~91&83~~97&97~~78&77~~84&84~~91&72~~78&71~~90&71~~94&72~~95&72~~92&72~~82&71~~92&71~~85\\
79~~88&80~~89&81~~90&83~~90&82~~98&79~~94&78~~81&85~~90&73~~84&73~~83&72~~80&73~~82&73~~95&73~~90&73~~79&72~~93\\
81~~86&82~~87&82~~89&84~~89&81~~99&85~~86&89~~91&86~~89&74~~80&74~~92&74~~83&74~~86&74~~96&74~~84&74~~94&73~~98\\
82~~85&83~~86&83~~88&85~~88&77~~80&81~~88&79~~99&87~~88&75~~88&75~~85&75~~89&69~~99&75~~90&75~~94&75~~87&75~~95\\
83~~84&84~~85&84~~87&86~~87&78~~79&80~~87&80~~96&51~~52&64~~71&64~~72&64~~73&64~~75&65~~66&65~~70&65~~72&65~~74\\
~~&~~&~~&~~&~~&~~&~~&~~&~~&~~&~~&~~&~~&~~&~~&~~\\
~~&~~&~~&~~&~~&~~&~~&~~&~table~24~&~&~&~&~&~&~~&\\
\end{tabular}

\end{sideways}
\end{center}}

\footnotesize{
\begin{center}
\begin{sideways}
\begin{tabular}{cccccccccccccccc}
34&35&36&37&38&39&40&41&42&43&44&45&46&47&48&49\\
77~~62&78~~62&76~~58&76~~59&76~~63&76~~61&76~~62&77~~57&76~~56&76~~57&77~~59&76~~55&77~~55&78~~51&76~~51&76~~53\\
80~~60&79~~56&77~~60&78~~56&78~~57&79~~63&81~~57&78~~60&77~~63&77~~58&78~~55&77~~53&79~~53&79~~55&78~~53&77~~51\\
81~~59&80~~63&78~~61&79~~60&79~~59&80~~58&82~~59&80~~56&82~~58&81~~56&81~~51&79~~51&80~~51&80~~53&79~~58&78~~58\\
83~~54&81~~61&79~~57&80~~61&81~~53&82~~57&83~~53&81~~58&83~~51&82~~51&82~~53&83~~60&85~~63&82~~56&80~~57&79~~62\\
85~~50&83~~55&80~~55&81~~62&84~~50&84~~54&84~~51&83~~62&84~~61&83~~59&83~~61&84~~58&86~~59&83~~57&82~~63&80~~59\\
86~~52&85~~52&82~~52&82~~55&88~~51&85~~53&85~~54&84~~53&85~~59&85~~55&88~~62&85~~57&88~~57&84~~59&84~~55&81~~55\\
87~~58&86~~50&85~~51&83~~58&90~~52&86~~51&86~~60&86~~55&86~~57&87~~53&89~~57&86~~62&90~~56&87~~54&87~~56&82~~61\\
88~~55&89~~51&87~~50&86~~53&91~~55&89~~50&88~~52&87~~51&87~~55&88~~60&90~~58&88~~59&92~~58&89~~60&90~~60&84~~60\\
90~~51&90~~53&89~~53&87~~52&93~~56&90~~55&89~~55&88~~54&88~~53&89~~52&91~~56&91~~63&93~~54&90~~61&91~~54&86~~56\\
92~~56&92~~54&90~~54&88~~50&94~~62&92~~52&91~~50&90~~50&93~~50&93~~63&92~~50&93~~52&94~~61&92~~62&92~~59&87~~63\\
93~~61&93~~60&92~~63&91~~51&95~~60&93~~59&92~~61&91~~52&94~~52&95~~50&94~~60&94~~50&95~~62&94~~63&93~~62&89~~54\\
95~~53&94~~57&95~~59&92~~57&97~~61&94~~56&94~~58&95~~61&95~~54&96~~54&96~~52&95~~56&97~~50&95~~58&95~~52&96~~57\\
96~~63&97~~58&96~~62&94~~54&98~~58&96~~60&98~~56&97~~63&97~~62&98~~61&97~~54&98~~54&98~~60&96~~50&96~~61&97~~52\\
98~~57&99~~59&99~~56&95~~63&99~~54&98~~62&99~~63&98~~59&99~~60&99~~62&98~~63&99~~61&99~~52&98~~52&98~~50&99~~50\\
64~~97&64~~87&64~~88&64~~93&64~~83&64~~95&64~~79&64~~94&64~~96&64~~84&64~~80&64~~82&64~~81&64~~86&64~~89&64~~92\\
66~~84&65~~84&65~~98&65~~89&65~~96&65~~81&65~~97&65~~79&65~~92&65~~94&65~~95&65~~87&65~~82&65~~88&65~~86&65~~91\\
67~~78&68~~77&67~~93&67~~97&67~~80&67~~87&66~~78&66~~82&66~~79&66~~90&66~~93&66~~97&66~~83&66~~91&66~~99&66~~94\\
68~~76&69~~96&68~~83&68~~90&68~~89&68~~99&69~~80&68~~96&68~~91&68~~78&68~~79&67~~92&67~~84&67~~81&67~~85&67~~90\\
69~~82&70~~82&70~~94&69~~77&69~~86&69~~78&70~~95&69~~89&69~~81&69~~92&69~~87&70~~90&69~~91&69~~85&68~~94&68~~85\\
70~~91&71~~95&71~~84&70~~98&70~~92&70~~97&71~~87&71~~93&70~~89&70~~80&71~~86&71~~89&70~~78&70~~99&69~~83&69~~93\\
71~~99&72~~98&72~~91&72~~99&71~~77&71~~91&72~~90&72~~76&71~~98&71~~97&74~~76&73~~80&71~~96&71~~76&71~~81&70~~88\\
72~~79&73~~76&73~~86&73~~96&72~~85&72~~77&73~~77&73~~85&73~~78&72~~86&73~~99&74~~81&72~~87&73~~93&72~~88&72~~83\\
73~~94&74~~88&74~~97&74~~85&74~~87&73~~88&74~~93&74~~99&74~~90&74~~91&72~~84&75~~78&73~~89&74~~77&73~~97&74~~95\\
74~~89&75~~91&75~~81&75~~84&75~~82&75~~83&75~~96&75~~92&75~~80&75~~79&70~~85&72~~96&75~~76&75~~97&75~~77&75~~98\\
65~~75&66~~67&66~~69&66~~71&66~~73&66~~74&67~~68&67~~70&67~~72&67~~73&67~~75&68~~69&68~~74&68~~72&70~~74&71~~73\\
~~&~~&~~&~~&~~&~~&~~&~~&~~&~~&~~&~~&~~&~~&~~&~~\\
~~&~~&~~&~~&~~&~~&~~&~~&~table~24~&~&~&~&~&~&~~&\\
\end{tabular}

\end{sideways}
\end{center}}

\footnotesize{
\begin{center}
\begin{sideways}

\begin{tabular}{ccccccccccccccccc}  \label{tbl1}
~~&~~&~~&~~&~~&~~&~~&~\textbf{BSTS($99$)}~&~~&~~&~~&~~&~~&~~&~~&~~&~~\\
1&2&3&4&5&6&7&8&9&10&11&12&13&14&15&16&17\\
50~~60&50~~69&56~~50&56~~52&56~~51&56~~53&56~~55&56~~54&56~~60&56~~61&60~~52&60~~53&60~~87&60~~90&60~~88&60~~93&60~~89\\
51~~61&51~~74&57~~51&57~~50&57~~54&57~~52&57~~53&57~~55&57~~69&57~~79&61~~53&61~~52&61~~97&61~~91&61~~87&61~~92&61~~99\\
52~~62&52~~70&58~~53&58~~54&58~~55&58~~50&58~~51&58~~52&58~~68&58~~62&62~~86&62~~87&62~~53&62~~89&62~~98&62~~95&62~~59\\
53~~79&53~~73&59~~52&59~~55&59~~53&59~~54&59~~50&59~~51&59~~70&59~~78&63~~85&63~~94&63~~52&63~~53&63~~99&63~~59&63~~58\\
54~~64&54~~72&54~~55&51~~53&50~~52&51~~55&52~~54&50~~53&50~~51&50~~54&64~~96&64~~95&64~~86&64~~52&64~~59&64~~58&64~~81\\
55~~77&55~~71&60~~63&60~~64&60~~65&60~~66&60~~67&60~~68&52~~55&53~~55&65~~97&65~~88&65~~98&65~~59&65~~58&65~~52&65~~54\\
65~~76&68~~75&62~~64&63~~65&64~~66&65~~67&66~~68&67~~69&53~~54&51~~52&66~~84&66~~93&66~~59&66~~58&66~~52&66~~55&66~~87\\
66~~75&67~~76&61~~65&62~~66&63~~67&64~~68&65~~69&66~~70&67~~71&60~~70&67~~98&67~~59&67~~58&67~~55&67~~53&67~~89&67~~52\\
67~~74&66~~77&66~~79&61~~67&62~~68&63~~69&64~~70&65~~71&66~~72&69~~71&68~~59&68~~51&68~~99&68~~94&68~~55&68~~80&68~~53\\
68~~73&65~~78&67~~78&68~~79&61~~69&62~~70&63~~71&64~~72&65~~73&68~~72&69~~82&69~~58&69~~85&69~~88&69~~54&69~~53&69~~55\\
69~~72&64~~79&68~~77&69~~78&70~~79&71~~61&72~~62&63~~73&64~~74&67~~73&70~~57&70~~50&70~~51&70~~96&70~~83&70~~54&70~~82\\
70~~71&61~~63&69~~76&70~~77&71~~78&72~~79&61~~73&62~~74&63~~75&66~~74&71~~56&71~~57&71~~50&71~~51&71~~81&71~~94&71~~88\\
63~~78&60~~62&70~~75&71~~76&72~~77&73~~78&74~~79&61~~75&62~~76&65~~75&72~~58&72~~56&72~~57&72~~95&72~~84&72~~91&72~~83\\
80~~81&80~~82&71~~74&72~~75&73~~76&74~~77&75~~78&76~~79&61~~77&64~~76&73~~50&73~~89&73~~56&73~~57&73~~51&73~~90&73~~84\\
82~~99&81~~83&72~~73&73~~74&74~~75&75~~76&76~~77&77~~78&78~~79&63~~77&74~~99&74~~92&74~~82&74~~56&74~~57&74~~96&74~~86\\
83~~98&84~~99&80~~89&80~~90&80~~83&80~~84&80~~85&80~~86&80~~87&80~~88&75~~83&75~~91&75~~81&75~~50&75~~56&75~~57&75~~85\\
84~~97&85~~98&88~~90&89~~91&82~~84&83~~85&84~~86&85~~87&86~~88&87~~89&76~~51&76~~90&76~~55&76~~54&76~~50&76~~56&76~~57\\
85~~96&86~~97&87~~91&88~~92&81~~85&82~~86&83~~87&84~~88&85~~89&86~~90&77~~81&77~~80&77~~54&77~~93&77~~86&77~~50&77~~56\\
86~~95&87~~96&86~~92&87~~93&86~~99&81~~87&82~~88&83~~89&84~~90&85~~91&78~~55&78~~54&78~~83&78~~92&78~~85&78~~51&78~~50\\
87~~94&88~~95&85~~93&86~~94&87~~98&88~~99&81~~89&82~~90&83~~91&84~~92&79~~54&79~~55&79~~84&79~~80&79~~82&79~~97&79~~51\\
88~~93&89~~94&84~~94&85~~95&88~~97&89~~98&90~~99&91~~81&82~~92&83~~93&80~~91&86~~96&80~~92&87~~97&80~~93&88~~98&80~~94\\
89~~92&90~~93&83~~95&84~~96&89~~96&90~~97&91~~98&92~~99&93~~81&82~~94&90~~92&85~~97&91~~93&86~~98&92~~94&87~~99&93~~95\\
90~~91&91~~92&82~~96&83~~97&90~~95&91~~96&92~~97&93~~98&94~~99&81~~95&89~~93&84~~98&90~~94&85~~99&91~~95&81~~86&92~~96\\
56~~57&56~~58&81~~97&82~~98&91~~94&92~~95&93~~96&94~~97&95~~98&96~~99&88~~94&83~~99&89~~95&84~~81&90~~96&82~~85&91~~97\\
58~~59&57~~59&98~~99&81~~99&92~~93&93~~94&94~~95&95~~96&96~~97&97~~98&87~~95&81~~82&88~~96&82~~83&89~~97&83~~84&90~~98\\
~~&~~&~~&~~&~~&~~&~~&~~&~~&~~&~~&~~&~~&~~&~~&~~&~~\\
~~&~~&~~&~~&~~&~~&~~&table~25&~~&~&~&~&~&~&~~&\\
\end{tabular}

\end{sideways}
\end{center}}

\footnotesize{
\begin{center}
\begin{sideways}
\begin{tabular}{cccccccccccccccc}
18&19&20&21&22&23&24&25&26&27&28&29&30&31&32&33\\
60~~59&60~~58&60~~57&60~~96&60~~55&60~~94&60~~84&60~~80&60~~85&60~~54&60~~51&80~~50&80~~51&80~~67&80~~70&80~~52\\
61~~58&61~~59&61~~50&61~~95&61~~83&61~~55&61~~93&61~~98&61~~54&61~~80&61~~57&81~~66&81~~72&81~~65&81~~60&81~~67\\
62~~93&62~~81&62~~56&62~~57&62~~51&62~~80&62~~55&62~~54&62~~94&62~~99&62~~50&82~~77&82~~50&82~~62&82~~51&82~~56\\
63~~95&63~~90&63~~51&63~~50&63~~57&63~~96&63~~54&63~~55&63~~86&63~~81&63~~56&83~~57&83~~60&83~~50&83~~74&83~~54\\
64~~53&64~~89&64~~80&64~~56&64~~50&64~~57&64~~51&64~~97&64~~55&64~~98&64~~82&84~~63&84~~57&84~~77&84~~50&84~~62\\
65~~80&65~~53&65~~95&65~~51&65~~56&65~~50&65~~57&65~~99&65~~93&65~~55&65~~91&85~~79&85~~73&85~~57&85~~71&85~~50\\
66~~94&66~~88&66~~53&66~~80&66~~54&66~~56&66~~50&66~~57&66~~51&66~~82&66~~83&86~~56&86~~58&86~~51&86~~59&86~~55\\
67~~96&67~~82&67~~94&67~~54&67~~92&67~~51&67~~56&67~~50&67~~57&67~~97&67~~90&87~~52&87~~67&87~~58&87~~55&87~~68\\
68~~52&68~~83&68~~96&68~~97&68~~84&68~~98&68~~85&68~~56&68~~50&68~~57&68~~54&88~~61&88~~52&88~~78&88~~57&88~~51\\
69~~97&69~~52&69~~93&69~~94&69~~91&69~~95&69~~92&69~~51&69~~56&69~~59&69~~84&89~~59&89~~68&89~~52&89~~75&89~~65\\
70~~55&70~~84&70~~97&70~~98&70~~85&70~~99&70~~86&70~~81&70~~58&70~~56&70~~53&90~~64&90~~69&90~~56&90~~68&90~~61\\
71~~92&71~~87&71~~52&71~~53&71~~90&71~~54&71~~91&71~~96&71~~59&71~~58&71~~89&91~~55&91~~74&91~~79&91~~52&91~~66\\
72~~51&72~~86&72~~55&72~~93&72~~53&72~~94&72~~90&72~~59&72~~87&72~~50&72~~52&92~~58&92~~55&92~~64&92~~76&92~~57\\
73~~91&73~~54&73~~98&73~~52&73~~86&73~~81&73~~59&73~~58&73~~92&73~~96&73~~55&93~~54&93~~71&93~~55&93~~56&93~~59\\
74~~90&74~~50&74~~54&74~~55&74~~52&74~~59&74~~58&74~~95&74~~88&74~~53&74~~85&94~~65&94~~59&94~~61&94~~58&94~~64\\
75~~54&75~~55&75~~92&75~~99&75~~59&75~~58&75~~52&75~~82&75~~53&75~~51&75~~88&95~~78&95~~70&95~~59&95~~73&95~~53\\
76~~98&76~~85&76~~99&76~~59&76~~58&76~~52&76~~87&76~~53&76~~91&76~~83&76~~86&96~~62&96~~56&96~~66&96~~53&96~~58\\
77~~57&77~~51&77~~91&77~~58&77~~87&77~~53&77~~88&77~~83&77~~89&77~~52&77~~59&97~~76&97~~75&97~~53&97~~54&97~~63\\
78~~56&78~~57&78~~58&78~~81&78~~89&78~~82&78~~53&78~~52&78~~90&78~~84&78~~87&98~~53&98~~54&98~~63&98~~69&98~~78\\
79~~50&79~~56&79~~59&79~~92&79~~88&79~~93&79~~89&79~~94&79~~52&79~~95&79~~58&99~~51&99~~53&99~~54&99~~72&99~~79\\
81~~88&80~~95&81~~90&82~~91&80~~97&83~~92&80~~98&84~~93&80~~99&85~~94&80~~96&60~~71&76~~66&60~~72&77~~67&60~~73\\
82~~87&94~~96&82~~89&83~~90&96~~98&84~~91&97~~99&85~~92&81~~98&86~~93&95~~97&70~~72&65~~77&71~~73&66~~78&72~~74\\
83~~86&93~~97&83~~88&84~~89&95~~99&85~~90&81~~96&86~~91&82~~97&87~~92&94~~98&69~~73&64~~78&70~~74&65~~79&71~~75\\
84~~85&92~~98&84~~87&85~~88&94~~81&86~~89&95~~82&87~~90&96~~83&88~~91&93~~99&68~~74&63~~79&69~~75&61~~64&70~~76\\
89~~99&91~~99&85~~86&86~~87&82~~93&87~~88&83~~94&88~~89&84~~95&89~~90&81~~92&67~~75&61~~62&68~~76&62~~63&69~~77\\
~~&~~&~~&~~&~~&~~&~~&~~&~~&~~&~~&~~&~~&~~&~~&~~\\
~~&~~&~~&~~&~~&~~&~~&~~&~table~25~&~&~&~&~&~&~~&\\
\end{tabular}

\end{sideways}
\end{center}}

\footnotesize{
\begin{center}
\begin{sideways}
\begin{tabular}{cccccccccccccccc}
34&35&36&37&38&39&40&41&42&43&44&45&46&47&48&49\\
80~~73&80~~56&80~~59&80~~71&80~~58&80~~63&80~~74&80~~53&80~~57&80~~54&80~~55&80~~69&80~~76&80~~78&80~~72&80~~75\\
81~~50&81~~61&81~~74&81~~68&81~~79&81~~56&81~~52&81~~69&81~~53&81~~55&81~~54&81~~58&81~~57&81~~59&81~~51&81~~76\\
82~~57&82~~68&82~~72&82~~54&82~~52&82~~71&82~~61&82~~58&82~~59&82~~73&82~~76&82~~53&82~~55&82~~63&82~~65&82~~60\\
83~~58&83~~64&83~~55&83~~59&83~~53&83~~69&83~~73&83~~67&83~~63&83~~65&83~~79&83~~62&83~~56&83~~51&83~~52&83~~71\\
84~~76&84~~51&84~~75&84~~64&84~~61&84~~58&84~~56&84~~71&84~~54&84~~67&84~~52&84~~65&84~~59&84~~55&84~~53&84~~74\\
85~~51&85~~65&85~~52&85~~56&85~~55&85~~54&85~~58&85~~59&85~~62&85~~66&85~~61&85~~64&85~~72&85~~53&85~~67&85~~77\\
86~~54&86~~50&86~~78&86~~67&86~~75&86~~68&86~~53&86~~65&86~~61&86~~71&86~~60&86~~52&86~~79&86~~66&86~~57&86~~69\\
87~~69&87~~57&87~~70&87~~50&87~~59&87~~65&87~~51&87~~64&87~~79&87~~53&87~~56&87~~63&87~~74&87~~54&87~~54&87~~73\\
88~~59&88~~58&88~~50&88~~53&88~~56&88~~67&88~~76&88~~70&88~~55&88~~72&88~~62&88~~68&88~~73&88~~58&88~~63&88~~64\\
89~~74&89~~66&89~~57&89~~51&89~~50&89~~70&89~~54&89~~55&89~~56&89~~69&89~~63&89~~61&89~~53&89~~77&89~~76&89~~72\\
90~~70&90~~62&90~~53&90~~65&90~~57&90~~50&90~~59&90~~54&90~~52&90~~51&90~~75&90~~66&90~~58&90~~67&90~~55&90~~79\\
91~~56&91~~59&91~~51&91~~70&91~~78&91~~53&91~~62&91~~68&91~~50&91~~57&91~~64&91~~54&91~~60&91~~56&91~~58&91~~63\\
92~~52&92~~63&92~~77&92~~62&92~~60&92~~66&92~~50&92~~72&92~~51&92~~68&92~~53&92~~70&92~~54&92~~52&92~~59&92~~65\\
93~~75&93~~53&93~~58&93~~63&93~~51&93~~64&93~~57&93~~73&93~~76&93~~70&93~~50&93~~67&93~~78&93~~76&93~~74&93~~68\\
94~~53&94~~52&94~~54&94~~57&94~~73&94~~51&94~~77&94~~50&94~~60&94~~74&94~~78&94~~55&94~~75&94~~57&94~~56&94~~70\\
95~~77&95~~54&95~~71&95~~58&95~~76&95~~55&95~~60&95~~66&95~~75&95~~50&95~~51&95~~56&95~~52&95~~65&95~~68&95~~67\\
96~~55&96~~69&96~~76&96~~52&96~~54&96~~72&96~~79&96~~57&96~~78&96~~59&96~~77&96~~51&96~~50&96~~73&96~~75&96~~61\\
97~~71&97~~55&97~~56&97~~66&97~~72&97~~52&97~~78&97~~51&97~~74&97~~58&97~~57&97~~59&97~~77&97~~50&97~~50&97~~62\\
98~~72&98~~79&98~~60&98~~55&98~~74&98~~59&98~~75&98~~56&98~~77&98~~52&98~~58&98~~57&98~~51&98~~64&98~~71&98~~66\\
99~~60&99~~67&99~~73&99~~69&99~~77&99~~57&99~~55&99~~52&99~~58&99~~56&99~~59&99~~50&99~~71&99~~75&99~~66&99~~78\\
78~~68&60~~74&69~~79&60~~76&62~~71&60~~77&63~~72&60~~78&64~~73&60~~79&65~~74&60~~75&61~~70&60~~61&60~~69&56~~59\\
67~~79&73~~75&68~~61&75~~77&70~~63&76~~78&71~~64&77~~79&72~~65&78~~61&73~~66&74~~76&69~~62&62~~79&61~~79&57~~58\\
66~~61&72~~76&62~~67&74~~78&64~~69&75~~79&65~~70&61~~76&66~~71&62~~77&67~~72&73~~77&63~~68&69~~74&62~~78&50~~55\\
62~~65&71~~77&63~~66&73~~79&65~~68&61~~74&66~~69&62~~75&67~~70&63~~76&68~~71&72~~78&64~~67&68~~70&70~~73&51~~54\\
63~~64&70~~78&64~~65&61~~72&66~~67&62~~73&67~~68&63~~74&68~~69&64~~75&69~~70&71~~79&65~~66&71~~72&64~~77&52~~53\\
~~&~~&~~&~~&~~&~~&~~&~~&~~&~~&~~&~~&~~&~~&~~&~~\\
~~&~~&~~&~~&~~&~~&~~&~~&~table~25~&~&~&~&~&~&~~&\\
\end{tabular}

\end{sideways}
\end{center}}

\footnotesize{
\begin{center}
\begin{sideways}

\begin{tabular}{ccccccccccccccccc}  \label{tbl1}
~~&~~&~~&~~&~~&~~&~~&~\textbf{BSTS($99$)}~&~~&~~&~~&~~&~~&~~&~~&~~&~~\\
1&2&3&4&5&6&7&8&9&10&11&12&13&14&15&16&17\\
50~~60&50~~61&60~~77&60~~67&60~~72&60~~71&60~~74&60~~75&60~~76&60~~69&60~~51&60~~53&60~~54&60~~55&60~~56&60~~57&60~~58\\
51~~61&51~~62&61~~64&61~~68&61~~69&61~~66&61~~72&61~~74&61~~70&61~~78&61~~80&61~~54&61~~98&61~~53&61~~55&61~~56&61~~57\\
52~~62&52~~60&62~~79&62~~74&62~~76&62~~78&62~~68&62~~65&62~~75&62~~77&62~~81&62~~50&62~~55&62~~54&62~~99&62~~97&62~~56\\
53~~63&53~~64&63~~71&63~~69&63~~75&63~~72&63~~79&63~~78&63~~74&63~~76&63~~98&63~~51&63~~52&63~~50&63~~54&63~~55&63~~80\\
54~~64&54~~67&65~~73&64~~79&64~~78&64~~73&64~~76&64~~70&64~~77&64~~75&64~~96&64~~52&64~~51&64~~56&64~~50&64~~59&64~~55\\
55~~65&55~~69&66~~74&65~~75&65~~70&65~~77&65~~69&66~~67&65~~72&65~~74&65~~95&65~~59&65~~53&65~~51&65~~52&65~~50&65~~54\\
56~~66&56~~65&67~~75&66~~78&66~~79&67~~76&66~~70&68~~76&66~~71&66~~73&66~~82&66~~81&66~~99&66~~52&66~~51&66~~53&66~~50\\
57~~67&57~~68&68~~72&70~~73&67~~73&68~~69&67~~77&69~~77&67~~79&67~~72&67~~59&67~~96&67~~97&67~~87&67~~86&67~~51&67~~52\\
58~~68&58~~63&69~~78&71~~76&68~~77&70~~74&71~~78&71~~72&68~~78&68~~71&68~~87&68~~92&68~~91&68~~84&68~~53&68~~52&68~~51\\
59~~69&59~~66&70~~76&72~~77&71~~74&75~~79&73~~75&73~~79&69~~73&70~~79&69~~58&69~~80&69~~50&69~~94&69~~97&69~~54&69~~98\\
70~~71&70~~75&50~~51&50~~52&50~~53&50~~54&50~~55&50~~56&50~~57&50~~58&70~~50&70~~56&70~~84&70~~86&70~~93&70~~96&70~~53\\
72~~75&71~~73&52~~59&51~~53&51~~55&51~~57&51~~59&51~~52&51~~54&51~~56&71~~93&71~~82&71~~83&71~~99&71~~81&71~~95&71~~94\\
73~~74&72~~74&53~~58&54~~59&52~~54&52~~56&52~~58&53~~59&52~~53&52~~55&72~~56&72~~84&72~~86&72~~59&72~~82&72~~80&72~~83\\
76~~79&76~~77&54~~57&55~~58&56~~59&53~~55&53~~57&54~~58&55~~59&53~~54&73~~55&73~~97&73~~96&73~~81&73~~94&73~~86&73~~82\\
77~~78&78~~79&55~~56&56~~57&57~~58&58~~59&54~~56&55~~57&56~~58&57~~59&74~~91&74~~87&74~~57&74~~82&74~~96&74~~58&74~~59\\
80~~82&80~~84&80~~85&80~~86&80~~88&80~~89&80~~90&80~~91&80~~92&80~~99&75~~54&75~~57&75~~92&75~~97&75~~58&75~~81&75~~95\\
81~~83&81~~93&81~~98&81~~99&81~~95&81~~94&81~~92&81~~85&81~~91&81~~90&76~~57&76~~88&76~~94&76~~58&76~~59&76~~85&76~~86\\
94~~99&82~~90&82~~94&82~~92&82~~99&82~~98&82~~86&82~~87&82~~93&82~~83&77~~53&77~~89&77~~58&77~~57&77~~87&77~~88&77~~90\\
84~~92&83~~97&83~~95&83~~88&83~~93&83~~92&83~~91&83~~99&83~~90&84~~89&78~~94&78~~58&78~~56&78~~85&78~~57&78~~87&78~~89\\
85~~97&94~~96&84~~90&84~~91&84~~96&84~~95&84~~87&84~~93&84~~94&85~~92&79~~52&79~~55&79~~59&79~~95&79~~89&79~~90&79~~87\\
86~~93&85~~89&86~~87&85~~93&85~~87&85~~90&85~~99&86~~97&85~~86&86~~88&83~~85&83~~98&80~~81&80~~98&80~~83&82~~84&81~~88\\
87~~98&86~~98&88~~93&87~~96&86~~90&86~~91&88~~97&88~~92&87~~99&87~~94&84~~97&85~~91&82~~95&83~~96&84~~88&83~~89&84~~85\\
95~~96&87~~91&89~~92&89~~98&89~~97&87~~93&89~~94&89~~95&88~~98&91~~95&86~~99&86~~94&85~~88&88~~91&85~~95&91~~94&91~~99\\
89~~91&88~~95&91~~97&90~~97&91~~98&88~~96&93~~96&90~~94&89~~96&93~~98&88~~89&90~~99&87~~90&89~~90&90~~98&92~~98&92~~96\\
88~~90&92~~99&96~~99&94~~95&92~~94&97~~99&95~~98&96~~98&95~~97&96~~97&90~~92&93~~95&89~~93&92~~93&91~~92&93~~99&93~~97\\
~~&~~&~~&~~&~~&~~&~~&~~&~~&~~&~~&~~&~~&~~&~~&~~&~~\\
~~&~~&~~&~~&~~&~~&~~&table~26&~~&~&~&~&~&~&~~&\\
\end{tabular}

\end{sideways}
\end{center}}

\footnotesize{
\begin{center}
\begin{sideways}
\begin{tabular}{cccccccccccccccc}
18&19&20&21&22&23&24&25&26&27&28&29&30&31&32&33\\
60~~59&60~~80&60~~98&60~~97&60~~83&60~~96&60~~85&60~~86&60~~92&60~~88&60~~89&80~~50&80~~52&80~~51&80~~53&80~~57\\
61~~97&61~~82&61~~95&61~~94&61~~89&61~~92&61~~88&61~~52&61~~96&61~~59&61~~58&81~~51&81~~50&81~~52&81~~60&81~~53\\
62~~83&62~~58&62~~85&62~~88&62~~98&62~~93&62~~89&62~~53&62~~59&62~~84&62~~57&82~~75&82~~51&82~~50&82~~52&82~~65\\
63~~82&63~~94&63~~59&63~~86&63~~88&63~~85&63~~93&63~~84&63~~57&63~~56&63~~96&83~~53&83~~70&83~~73&83~~50&83~~51\\
64~~80&64~~89&64~~92&64~~58&64~~94&64~~91&64~~84&64~~95&64~~99&64~~57&64~~93&84~~52&84~~53&84~~71&84~~51&84~~50\\
65~~88&65~~84&65~~99&65~~98&65~~58&65~~83&65~~57&65~~90&65~~85&65~~89&65~~92&85~~67&85~~68&85~~53&85~~71&85~~72\\
66~~54&66~~55&66~~97&66~~57&66~~84&66~~94&66~~58&66~~83&66~~88&66~~92&66~~87&86~~78&86~~77&86~~54&86~~79&86~~58\\
67~~50&67~~53&67~~55&67~~56&67~~99&67~~58&67~~81&67~~89&67~~84&67~~93&67~~91&87~~64&87~~65&87~~69&87~~58&87~~70\\
68~~55&68~~50&68~~54&68~~59&68~~56&68~~90&68~~83&68~~88&68~~97&68~~81&68~~99&88~~70&88~~78&88~~58&88~~67&88~~71\\
69~~51&69~~52&69~~83&69~~96&69~~57&69~~56&69~~90&69~~82&69~~91&69~~53&69~~85&89~~63&89~~58&89~~76&89~~66&89~~68\\
70~~52&70~~51&70~~57&70~~80&70~~54&70~~55&70~~59&70~~58&70~~98&70~~82&70~~95&90~~58&90~~60&90~~61&90~~62&90~~59\\
71~~53&71~~54&71~~51&71~~52&71~~50&71~~57&71~~55&71~~56&71~~58&71~~87&71~~59&91~~61&91~~66&91~~65&91~~59&91~~56\\
72~~95&72~~57&72~~52&72~~51&72~~53&72~~50&72~~54&72~~98&72~~90&72~~58&72~~55&92~~69&92~~63&92~~59&92~~56&92~~67\\
73~~57&73~~59&73~~58&73~~53&73~~51&73~~52&73~~50&73~~54&73~~56&73~~99&73~~90&93~~54&93~~59&93~~56&93~~69&93~~79\\
74~~86&74~~93&74~~53&74~~83&74~~52&74~~51&74~~80&74~~50&74~~54&74~~55&74~~56&94~~59&94~~56&94~~57&94~~77&94~~75\\
75~~87&75~~86&75~~56&75~~55&75~~90&75~~53&75~~51&75~~59&75~~50&75~~98&75~~52&95~~56&95~~55&95~~68&95~~73&95~~69\\
76~~56&76~~95&76~~90&76~~91&76~~81&76~~54&76~~52&76~~51&76~~55&76~~50&76~~53&96~~62&96~~71&96~~66&96~~55&96~~54\\
77~~85&77~~56&77~~93&77~~54&77~~59&77~~80&77~~91&77~~55&77~~51&77~~52&77~~50&97~~57&97~~79&97~~78&97~~54&97~~55\\
78~~92&78~~91&78~~50&78~~93&78~~55&78~~59&78~~53&78~~99&78~~52&78~~51&78~~54&98~~71&98~~57&98~~55&98~~74&98~~73\\
79~~58&79~~92&79~~91&79~~50&79~~85&79~~81&79~~56&79~~57&79~~53&79~~54&79~~51&99~~55&99~~54&99~~72&99~~57&99~~52\\
81~~89&83~~87&80~~96&81~~84&82~~91&82~~97&82~~96&80~~87&80~~95&80~~97&80~~94&60~~65&61~~62&60~~62&61~~63&60~~61\\
84~~99&85~~98&81~~82&82~~85&80~~93&84~~98&86~~92&81~~96&81~~87&83~~94&81~~86&66~~68&64~~67&63~~64&64~~65&62~~63\\
90~~93&88~~99&84~~86&87~~92&86~~96&86~~89&87~~95&85~~94&82~~89&85~~96&82~~88&72~~73&69~~72&67~~70&68~~70&64~~66\\
91~~96&81~~97&87~~89&89~~99&87~~97&87~~88&94~~97&91~~93&83~~86&86~~95&83~~84&74~~76&73~~76&74~~79&75~~78&74~~77\\
94~~98&90~~96&88~~94&90~~95&92~~95&95~~99&98~~99&92~~97&93~~94&90~~91&97~~98&77~~79&74~~75&75~~77&72~~76&76~~78\\
~~&~~&~~&~~&~~&~~&~~&~~&~~&~~&~~&~~&~~&~~&~~&~~\\
~~&~~&~~&~~&~~&~~&~~&~~&~table~26~&~&~&~&~&~&~~&\\
\end{tabular}
\end{sideways}
\end{center}}

\footnotesize{
\begin{center}
\begin{sideways}
\begin{tabular}{cccccccccccccccc}
34&35&36&37&38&39&40&41&42&43&44&45&46&47&48&49\\
80~~58&80~~59&80~~66&80~~65&80~~67&80~~68&80~~54&80~~71&80~~73&80~~75&80~~78&80~~76&80~~56&80~~62&80~~55&80~~79\\
81~~57&81~~58&81~~59&81~~61&81~~63&81~~64&81~~65&81~~69&81~~70&81~~72&81~~74&81~~56&81~~77&81~~55&81~~54&81~~78\\
82~~53&82~~57&82~~58&82~~77&82~~78&82~~79&82~~59&82~~60&82~~62&82~~64&82~~76&82~~55&82~~67&82~~54&82~~56&82~~68\\
83~~52&83~~75&83~~77&83~~76&83~~79&83~~78&83~~61&83~~59&83~~57&83~~63&83~~55&83~~64&83~~54&83~~56&83~~58&83~~67\\
84~~73&84~~74&84~~75&84~~78&84~~76&84~~60&84~~57&84~~79&84~~59&84~~61&84~~56&84~~54&84~~55&84~~58&84~~77&84~~69\\
85~~50&85~~51&85~~52&85~~73&85~~74&85~~75&85~~64&85~~57&85~~56&85~~55&85~~54&85~~59&85~~58&85~~61&85~~66&85~~70\\
86~~51&86~~50&86~~53&86~~52&86~~61&86~~57&86~~62&86~~64&86~~55&86~~56&86~~65&86~~66&86~~59&86~~68&86~~69&86~~71\\
87~~72&87~~52&87~~50&87~~51&87~~53&87~~73&87~~55&87~~56&87~~54&87~~60&87~~61&87~~62&87~~63&87~~59&87~~57&87~~76\\
88~~74&88~~53&88~~55&88~~50&88~~51&88~~52&88~~56&88~~54&88~~72&88~~73&88~~75&88~~69&88~~79&88~~57&88~~59&88~~64\\
89~~59&89~~56&89~~51&89~~53&89~~50&89~~54&89~~52&89~~55&89~~69&89~~70&89~~71&89~~72&89~~57&89~~73&89~~74&89~~75\\
90~~56&90~~63&90~~64&90~~55&90~~54&90~~50&90~~51&90~~52&90~~53&90~~78&90~~67&90~~57&90~~66&90~~70&90~~71&90~~74\\
91~~75&91~~73&91~~71&91~~54&91~~55&91~~51&91~~50&91~~53&91~~52&91~~62&91~~57&91~~58&91~~60&91~~63&91~~70&91~~72\\
92~~70&92~~71&92~~54&92~~72&92~~57&92~~55&92~~53&92~~50&92~~51&92~~52&92~~58&92~~73&92~~74&92~~76&92~~62&92~~77\\
93~~76&93~~55&93~~72&93~~57&93~~75&93~~53&93~~60&93~~58&93~~50&93~~51&93~~52&93~~61&93~~65&93~~66&93~~68&93~~73\\
94~~55&94~~54&94~~74&94~~68&94~~52&94~~70&94~~67&94~~51&94~~58&94~~50&94~~53&94~~79&94~~62&94~~65&94~~72&94~~60\\
95~~54&95~~78&95~~57&95~~74&95~~77&95~~59&95~~58&95~~63&95~~66&95~~53&95~~50&95~~51&95~~52&95~~60&95~~67&95~~62\\
96~~68&96~~72&96~~56&96~~75&96~~59&96~~58&96~~79&96~~78&96~~76&96~~57&96~~77&96~~50&96~~51&96~~52&96~~53&96~~65\\
97~~77&97~~76&97~~70&97~~56&97~~58&97~~72&97~~74&97~~65&97~~71&97~~59&97~~64&97~~53&97~~50&97~~51&97~~52&97~~63\\
98~~78&98~~77&98~~79&98~~58&98~~56&98~~76&98~~68&98~~67&98~~64&98~~54&98~~59&98~~52&98~~53&98~~50&98~~51&98~~66\\
99~~63&99~~79&99~~69&99~~59&99~~70&99~~56&99~~76&99~~74&99~~60&99~~58&99~~51&99~~77&99~~75&99~~53&99~~50&99~~61\\
60~~66&60~~68&60~~78&60~~79&60~~73&61~~77&63~~70&61~~73&61~~75&65~~76&60~~70&60~~63&61~~71&64~~74&60~~64&50~~59\\
62~~64&65~~70&62~~67&62~~66&62~~71&62~~69&66~~75&62~~72&63~~77&66~~77&63~~68&65~~67&64~~68&67~~78&61~~79&51~~58\\
61~~65&62~~66&63~~65&63~~67&64~~72&63~~66&69~~71&70~~77&65~~79&67~~71&66~~72&68~~74&69~~76&69~~75&63~~73&52~~57\\
67~~69&61~~67&68~~73&64~~71&65~~68&65~~71&72~~78&66~~76&67~~68&68~~79&62~~73&70~~78&70~~72&71~~77&65~~78&53~~56\\
71~~79&64~~69&61~~76&69~~70&66~~69&67~~74&73~~77&68~~75&74~~78&69~~74&69~~79&71~~75&73~~78&72~~79&75~~76&54~~55\\
~~&~~&~~&~~&~~&~~&~~&~~&~~&~~&~~&~~&~~&~~&~~&~~\\
~~&~~&~~&~~&~~&~~&~~&~~&~table~26~&~&~&~&~&~&~~&\\
\end{tabular}
\end{sideways}
\end{center}}

\end{document}